\RequirePackage[l2tabu, orthodox]{nag}
\documentclass[a4paper,11pt,reqno]{amsart}

\usepackage{amssymb,amscd,amsthm}
\usepackage{mathrsfs}

\usepackage{url}
\usepackage{boondox-calo}
\usepackage{enumerate}

\usepackage[ansinew]{inputenc}
\usepackage[centertags]{amsmath}
\usepackage{epsfig,graphicx}
\usepackage{verbatim}

\usepackage[colorlinks,
           linkcolor=red,
           anchorcolor=green,
           citecolor=blue
           ]{hyperref}

\providecommand{\acknowledgement}
    {{\noindent{\textbf{{Acknowledgements:}}\quad{
    }}}}

\newcommand{\Z}{\mbox{$\mathbb{Z}$}}

\newcommand{\R}{\mbox{$\mathbb{R}$}}
\newcommand{\T}{\mbox{$\mathbb{T}$}}

   \def\cQ{\mathcal{Q}}

\newcommand{\bS}{\mathbb{S}}
\newcommand{\kS}{\mathfrak{S}}
\newcommand{\F}{\mathfrak{F}}
\newcommand{\kL}{\mathfrak{L}}
\newcommand{\E}{\mathcal{E}}
\newcommand{\FF}{\mathcal{F}}

\newcommand{\co}{\mathscr{O}}

\newenvironment{proof1}{\emph{The Proof of Proposition \ref{absolute}} :}{\hfill $\square$\par}

\newenvironment{proof2}{\emph{The Proof of item (ii) of Theorem \ref{final}} :}{\hfill $\square$\par}

\newenvironment{proof3}{\emph{The Proof of Proposition \ref{bound}} :}{\hfill $\square$\par}

\theoremstyle{plain}
\newtheorem{maintheorem}{Theorem}

\newtheorem*{teo*}{Theorem}
\newtheorem*{prop*}{Proposition}
\newtheorem*{cor*}{Corollary}
\newtheorem*{defi*}{Definition}
\newtheorem*{goal*}{Goal}

\newtheorem*{teoA'}{Theorem A'}

\newtheorem{teo}{Theorem}[section]
\newtheorem{thm}[teo]{Theorem}

\newtheorem{af}{Claim}
\newtheorem{lem}[teo]{Lemma}
\newtheorem{prop}[teo]{Proposition}

\newcommand{\bi}{\begin{itemize}}
\newcommand{\ei}{\end{itemize}}

\theoremstyle{definition}
\newtheorem{defi}[teo]{Definition}

\theoremstyle{remark}

\newtheorem{rmk}[teo]{Remark}

\numberwithin{equation}{section}

\author[X. Zhang]{Xiang Zhang}
\address{School of Mathematical Sciences, Peking University, Beijing 100871, People's Republic of China}
\email{julienhsiang@pku.edu.cn}


\title[On codimension one partially hyperbolic diffeomorphisms]
{On codimension one partially hyperbolic diffeomorphisms}

\begin{document}

\maketitle


\begin{abstract}
	We show that every codimension one partially hyperbolic diffeomorphism must support on $\T^{n}$. It is locally uniquely integrable and derived from a linear codimension one Anosov diffeomorphism. Moreover, this system is intrinsically ergodic, and the A. Katok's conjecture about the existence of ergodic measures with intermediate entropies holds for it.

\end{abstract}

\medskip

\section{Introduction}

Let $M$ be a $C^{\infty}$ closed Riemannian  $n$-manifold. A $C^{1}$-diffeomorphism $g: M \rightarrow M$ is an Anosov  diffeomorphism, if it admits a continuous $Tg$-invariant splitting $TM=E^{s} \oplus E^{u}$ such that for any $x\in M$, unit vectors $v^{s}\in E^{s}_{x}$ and  $v^{u}\in E^{u}_{x}$ we have
\begin{center}
$
  \parallel Tg(v^{s})\parallel < \tau < 1 < \tau^{-1} < \parallel Tg(v^{u})\parallel
$
\end{center}
for some constant $\tau<1$. When $g$ is codimension one, $i.e.$,  $dim(E^{s})=1$ or $dim(E^{u})=1$,  J. Franks and S. Newhouse proved that $g$ is topologically conjugate to a hyperbolic toral automorphism in  \cite{Fr} and \cite{Ne}. Then a natural question is, whether the similar result holds if the one dimensional sub-boundle has an intermediate behavior and the other $(n-1)$-dimensional sub-bundle is uniformly expanding or contracting. Our main issue will be the above question about the codimension one partially hyperbolic diffeomorphisms.

\begin{defi}\label{copa}
 A $C^{1}$-diffeomorphism $f: M \rightarrow M$ is called a codimension one  partially hyperbolic diffeomorphism, if $f$ or $f^{-1}$  admits a continuous $Tf$-invariant splitting $TM=E^{c} \oplus E^{u}$ and $E^{c}$ is one dimensional, and a function $\xi:M \rightarrow (1,+\infty)$ such that
 \begin{center}
$
  \parallel Tf(v^{c})\parallel<\xi(x)<\parallel Tf(v^{u})\parallel
$
\end{center}
for all $x\in M$ and unit vectors $v^{c}\in E^{c}_{x}$ and  $v^{u}\in E^{u}_{x}$.
\end{defi}

As is well known, S. Smale has shown how to construct a ``derived from Anosov'' diffeomorphism from an Anosov diffeomorphism of $\T^{2}$ \cite{Sma}, and the resulting diffeomorphism is structurally stable and has a \emph{DA}-attractor. Since the construction preserves the original stable foliation, we can make this \emph{DA}-diffeomorphism to be partially hyperbolic diffeomorphism. We refer to this construction here as the \emph{DA-construction}. See \cite{RP2} for the  two-dimensional case and more details.

In this paper, we show that every codimension one partially hyperbolic diffeomorphism is derived from a linear codimension one Anosov diffeomorphism  of a torus.

\begin{maintheorem}\label{final}

Let $N$ be a closed Riemannian manifold, and $f:N\rightarrow N$ be a codimension one partially hyperbolic diffeomorphism. Then,
\begin{enumerate}[(i).]
		\item $N$ is homeomorphic to $\T^{n}$;

        \item The distribution $E^{c}$ is locally uniquely integrable, $i.e.$ there is a foliation tangent to $E^{c}$ at each point and any piecewise $C^{1}$-curve tangent to $E^{c}$ must lie in a  unique leaf of the integral foliation (center foliation);
        \item $f$ is semiconjugated to a linear codimension one Anosov diffeomorphism $f_{\ast}$, where $f_{\ast}\in GL(n,\mathbb{Z})$ is the induced linear  transformation of $f$ on $\pi_{1}(N)$.
\end{enumerate}

\end{maintheorem}

\begin{rmk}\label{rmk0}
(1). In the next section, we will demonstrate that the torus $\mathbb{T}^{n}$ is the only closed smooth manifold which admits a $C^{0}$ \emph{hyperplanar foliation}, a codimension one foliation by hyperplanes. This can be considered as a generalization of Rosenberg's result \cite{Ha}, see Proposition \ref{torus}. As an immediate consequence we have that any codimension one partially hyperbolic diffeomorphism will only support on  $\T^{n}$. This result can be seen as an extension of the famous work of J. Franks and S. Newhouse on codimension one Anosov diffeomorphisms.

(2). From \cite{HPS}, there is always a unique foliation $\F^{u}$ tangent to $E^{u}$, called \emph{unstable foliation}. We will prove that any $C^{0,1}$ hyperplanar foliation is quasi-isometric (see Proposition \ref{iso1}), which will give another more succinct proof of the locally unique integrability of codimension one absolutely partially hyperbolic diffeomorphisms, see Proposition \ref{absolute} and its proof. A codimension one partially hyperbolic diffeomorphism  is \emph{absolute} if the function $\xi$ in Definition \ref{copa}  can be taken to be a constant on the whole $N$.

(3). In the following, we will know that the semiconjugacy only collapses center arcs, it maps the unstable foliation of $f$ into the unstable foliation of $f_{\ast}$ and maps the central foliation of $f$ into the stable foliation of $f_{\ast}$, see Proposition \ref{prei}. Therefore, combining Theorem \ref{final}, one can see that any codimension one partially hyperbolic  diffeomorphism can be obtained by imposing some $DA$-constructions on a linear codimension one Anosov diffeomorphism of a torus.

(4). In \cite{PS009}, E. Pujals and M. Sambarino showed that if $f$ is a $C^{2}$-diffeomorphism, $N$ admits a contractive codimension one dominated splitting and all the hyperbolic periodic points are of saddle type, then $N\cong\T^{n}$ and $f$ is an Anosov diffeomorphism. See \cite{PS009} for more details.

\end{rmk}

In the light of above theorem, we may give some dynamical characterizations about codimension one partially hyperbolic diffeomorphisms. We define a relation $\vartriangleleft$ on $\T^{n}$, for any $x,y\in \T^{n}$, $x\vartriangleleft y$ means that  for every $\epsilon>0$ there exists a set of points $z_{0}=x,\ldots, z_{n}=y$ such that $n\geq 1$ and $d(f(z_{i}),z_{i+1})<\epsilon$, $i\in\{0,1,\ldots,n-1\}$. Then the relation ``$\thickapprox$'' ($x \thickapprox y$ if and only if $x\vartriangleleft y$ and $y\vartriangleleft x$) in $\mbox{CR}(f)\triangleq\{x\in \T^{n}|\ x\vartriangleleft x\}$, the chain-recurrent set of $f$,  is an equivalence relation. Its equivalence classes are called  \emph{chain-recurrence classes}. And we call a compact invariant set $\cQ$ a \emph{quasi-attractor} for $f$ if it is a chain-recurrence class and there exists a decreasing sequence of open neighborhoods $\{V_{n}\}$ such that $\bigcap_{n} V_{n}=\cQ$ and $f(\overline{V_{n}})\subset V_{n}$. There always exists a quasi-attractor for $f$ due to the fundamental theorem of dynamical systems, Conley's theorem in \cite{Con}.

\begin{maintheorem}\label{dynachara}

Let $f:\T^{n}\rightarrow \T^{n}$ be a codimension one partially hyperbolic diffeomorphism, then there is a unique quasi-attractor $\mathcal{Q}$ of $f$. Furthermore, every chain-recurrence class different from $\mathcal{Q}$ is contained in a periodic center interval.

\end{maintheorem}

\begin{rmk}\label{rmk1}
The two dimensional case of the above theorem is given in Theorem 4.A.3. in \cite{RP2}.
\end{rmk}

It is worth mentioning that by Theorem \ref{final} and the proof of Theorem \ref{dynachara} we can conclude that $f$ is  \emph{intrinsically ergodic}, that is, $f$ has an unique entropy maximizing measure, which is an invariant measure such that its entropy is equal to the topological entropy of the system. In addition, we also show that any codimension one partially hyperbolic diffeomorphism admits ergodic measures whose entropies can be any value between zero and its topological entropy. In other words, the A. Katok's conjecture about the existence of ergodic measures with intermediate entropies is true for our setting.

\begin{maintheorem}\label{cord}

Let $f:\T^{n}\rightarrow \T^{n}$ be a codimension one partially hyperbolic diffeomorphism. Then

\begin{enumerate}[(i).]
\item $f$ is intrinsically ergodic. 
    Moreover, $(f,\mu)$ and $(f_{\ast},m)$ are isomorphic, where $m$ is the Lebesgue measure on $\T^{n}$;

\item For any $\lambda\in [0,\E_{\text {top}}(f)]$, there is an ergodic measure $\nu$ such that $\E_{\nu}(f)=\lambda$, where $\E_{\text {top}}(f)$ is the topological entropy of $f$ and $\E_{\nu}(f)$ is the measure entropy of $\nu$.

\end{enumerate}

\end{maintheorem}

The structure of this paper is as follows. In Section \ref{sec2} we prepare some results concerning geometric theory of foliations without holonomy. In subsection \ref{sec2.1},  there contains the existence of the transverse foliations and the global product structure, which are the key tools to prove the theorems above.  In subsection \ref{sec2.2}, we introduce some important properties of hyperplanar foliations, among which the main results are Proposition \ref{torus} and Proposition \ref{bound}, and item $(i)$ of Theorem \ref{final} follows from Proposition \ref{torus}. These propositions play an important role in proving the locally unique integrability of central distribution. Section  \ref{sec3} and Section \ref{sec4} are devoted to the proof of item $(ii)$ and item $(iii)$ of Theorem \ref{final},  respectively. In Section \ref{sec5}, we provide some specific dynamical characterizations about codimension one partially hyperbolic diffeomorphisms, the proof of Theorem \ref{dynachara} is shown there. At the end, in Section \ref{sec6} we give a proof of Theorem \ref{cord}.

\bigskip

\section{Geometric theory of codimension one foliations }\label{sec2}

For partially hyperbolic systems, there are always unique foliations $\F^{u}$ and $\F^{s}$ tangent to $E^{u}$ and $E^{s}$, whose leaves are homeomorphic to $\R^{dim(E^u)}$ and $\R^{dim(E^s)}$ \cite{HPS}, called (\emph{strong})  \emph{stable} and \emph{unstable foliations}, respectively. So these leaves of  $\F^{u}$ and $\F^{s}$ all have trivial holonomy group. And we call this kind of foliation is \emph{without holonomy}.

Therefore, in view of the partially hyperbolic systems considered in this paper are codimension one,  in this section we will introduce some preliminary knowledge about codimension one foliations without holonomy and show some important properties to facilitate subsequent proofs, in which we will focus on hyperplanar foliations.

\subsection{Preliminaries and Analysis}\label{sec2.1}

\begin{defi}
Let $M(M^{n})$ be a smooth closed  $n$-manifold. A \emph{codimension one foliation} $\kS$, is a decomposition of $M$ into a disjoint union of connect hypersurfaces, called the \emph{leaves} of $\kS$, together with a collection of charts $U_{i}$ covering $M$, with $\varphi_{i}:\R^{n-1}\times\R\rightarrow U_{i}$ a homeomorphism, such that the preimage of each component of a leaf intersected with $U_{i}$ is a horizontal hyperplane. In particular, $\kS$ is said to be \emph{hyperplanar} if all its leaves are homeomorphic to $\R^{n-1}$.

The foliation $\kS$ is $C^{r},\ 0 \leq r\leq \infty$, if the charts $(U_{i},\varphi_{i})$ can be chosen such that each $\varphi_{i}$ is a $C^{r}$-diffeomorphism. Besides, we say $\kS$ is $C^{0,r},\ 0 \leq r\leq \infty$, if the charts $(U_{i},\varphi_{i})$ can be chosen such that the restriction of each $\varphi_{i}$ to a horizontal hyperplane is a $C^{r}$-immersion and the tangent planes of leaves vary continuously. See \cite{CC} and \cite{HH1} for the general definition and more details of the geometric theory of foliations.

\end{defi}

In this subsection, the foliations we consider are $C^{0}$, if not noted explicitly. Although the $C^{0,1}$ case is enough for us to prove the main theorems in the introduction above. It is necessary to indicate that, in general, the regularity of $\F^{u}$ of $f$ is H\"{o}lder continuous. It can not be $C^{1}$ or have a better regularity, although $f$ is $C^{2}$ even smooth, see \cite{BP}, \cite{BDV} and references therein. However, when $f$ is a codimension one $C^{2}$ partially hyperbolic diffeomorphism,  $\F^{u}$ can be $C^{1}$, see \cite{PSW} and \cite{Pe}.

Let $\pi:\tilde{M}\rightarrow M$ be the universal covering map, throughout the paper we use $\tilde{A}$  to represent the $A$ lifted on the universal covering space. And denote the foliation transverse to $\kS$ as $\kS^{\pitchfork}$. $\kS(x)$ stands for the leaf of $\kS$ that contains $x$. The transverse foliation is one of the main tools in subsequent proofs, so  we need to discuss its existence in detail. First, we give the following definition of transversality in $C^r$ sense, where $r\geq 0$.

\begin{defi}[Transversality] Let $\kS$ be a codimension $k$ foliation on a closed manifold $\Lambda^{m}$. A submanifold $\Lambda'$ of $\Lambda$ is called transverse to $\kS$ if the inclusion $i: \Lambda'\rightarrow \Lambda$ is transverse to $\kS$, $i.e.$, for any $x'\in \Lambda'$ if there is a neighborhood $U'$ of $x'$ and a distinguished map $(U,\phi)$ of $\kS$, $i(U')\subset U$, such that $\phi\circ i:U'\rightarrow\R^{k}$ is a submersion. In particular, a foliation $\kS'$ on $\Lambda$ is said to be transverse to $\kS$ if each leaf of $\kS'$ is transverse to $\kS$.

In the $C^{0}$ case, a submersion $\varphi:\Lambda^{m}\rightarrow M^{n}$, $n\leq m$, is locally of the form $\varphi=P\circ h$, where $h$ is a homeomorphism in $\R^{m}$ and $P:\R^{m}=\R^{m-n}\times \R^{n}\rightarrow \R^{n}$ is the canonical projection.

\end{defi}

In general, a foliation does not admit any transverse foliation of complementary dimension, see section 2.3 in Chapter II of \cite{HH1}. While in the case of codimension one,  there is an optimistic result, even $\kS$ is $C^{0}$, the existence of $\kS^{\pitchfork}$ is guaranteed by the following theorem.

\begin{thm}[\cite{HH2}, Theorem IV.1.1.2]\label{existrans}
Consider a $C^{r}\ (0\leq r\leq \infty)$ codimension one foliation $(M,\kS)$, possibly with boundary, and let $K\subset M$ be a compact set.  Suppose that a transverse foliation  $\kS_{K}^{\pitchfork}$ of $\kS$ is given on some open neighborhood of $K$. Then there exists a transverse foliation  $\kS^{\pitchfork}$ of $\kS$ such that $\kS_{K}^{\pitchfork}$ and $\kS^{\pitchfork}$ agree on some open neighborhood of $K$.
In particular, $\kS^{\pitchfork}$ has a certain set of closed transversals of $\kS$ as leaves.
\end{thm}

Indeed, when $\kS$ is $C^{0,1}$, the existence of $\kS^{\pitchfork}$ is easy to check. When $\kS$ is transversally orientable, we can choose  a $C^{1}$-vector field without any singularity transverse to $T\kS$ in $TM$ (by slightly disturbing the orthogonal subbundle to $T\kS$), it will integrate to a foliation transverse to $\kS$. If $\kS$ is not transversally orientable, consider a double cover, then just take a $C^{1}$-line field transverse to $T\kS$ and invariant under deck transformations.

\begin{defi}
  Given two transverse foliations $\F_{1}$ and $\F^{\pitchfork}_{1}$ on a manifold $M$, we say they have  \emph{global product structure} if for any two points $x,\ y\in \tilde{M}$, the leaves $\tilde{\F_{1}}(x)$ and $\tilde{\F}^{\pitchfork}_{1}(y)$ intersect in a unique point.
\end{defi}

When $\F_{1}$ is a $C^{r},\ r\geq 0$ foliation, we have the following theorem,

\begin{thm}[\cite{HH2}, Theorem VIII.2.2.1]\label{gps}
Consider a codimension one foliation $\F_{1}$ of a closed manifold $M$ such that all the leaves of $\F_{1}$ have trival holonomy group. Then, for every foliation $\F^{\pitchfork}_{1}$ transverse to $\F_{1}$, $\F_{1}$ and $\F^{\pitchfork}_{1}$ have global product structure. Furthermore, $\tilde{M}$ is homeomorphic to $\tilde{F}\times \R$, where $\tilde{F}$ is the universal covering of a leaf $F$ of $\F_{1}$.
\end{thm}

Actually, if $\F_{1}$ has better regularity, we can use other methods to prove Theorem \ref{gps}. For example, when $\F_{1}$ is $C^{r},\ r\geq 2$, we can apply  Theorem $9.2.1$ in book \cite{CC} to construct a $C^{0}$ leaf-preserving flow, which implies that $\F_{1}$ and $\F_{1}^{\pitchfork}$ have  global product structure. And then $\tilde{M}$ is even  diffeomorphic to $\tilde{F}\times \R$, not only homeomorphic, see \cite{No}\cite{HI}.

The global product structure often implies a lot of geometric information about the ambient manifold and topological characterization of foliations  themselves. We give the following several typical propositions that can reflect this fact, all of which can be regarded as corollaries of Theorem \ref{gps}.

\begin{prop}\label{glho}
Let $\F_{1}$ and $\F^{\pitchfork}_{1}$ be as in Theorem \ref{gps}, then $\tilde{M}$ is homeomorphic to the product of $\tilde{\F}_{1}(x) \times \tilde{\F}^{\pitchfork}_{1}(x)$ for any $x \in \tilde{M}$. In particular, all the leaves of $\tilde{\F}_{1}(\tilde{\F}^{\pitchfork}_{1})$  must be simply connected and they are homeomorphic to each other.

\begin{proof}
   Construct the map $\phi:\tilde{\F}_{1}(x) \times \tilde{\F}^{\pitchfork}_{1}(x)\rightarrow \tilde{M}$, $\phi(y,z)=\tilde{\F}^{\pitchfork}_{1}(y) \cap \tilde{\F}_{1}(z)\in \tilde{M}$.
By the global product structure it is easy to know that $\phi$ is well defined and bijective. The map is also continuous because of the continuity of foliations. By the invariance of domain theorem, $\phi^{-1}$ is continuous and $\phi$ is a global homeomorphism as desired.

  \end{proof}
\end{prop}

\begin{prop}\label{leafspace}
Let $\F_{1}$ and $\F^{\pitchfork}_{1}$  be as in Theorem \ref{gps}, the leaf space $\FF_{1}\triangleq \tilde{M}/\tilde{\F}_{1}$ $($ $i.e.$, $\FF_{1}$ is the quotient of $\tilde{M}$ by the equivalence relation $\thicksim$: $x\thicksim y$ if and only if $x,\ y$ are in the same leaf of $\tilde{\F}_{1}$ for any $x,\ y\in\tilde{M}$ $)$ of $\tilde{\F}_{1}$ is $\R$.
\begin{proof}
The leaf space of $\tilde{\F}_{1}$ with the quotient topology has the structure of a (possibly non-Hausdorff) one-dimensional ``manifold''($i.e.$, a space locally homeomorphic to $\R$ and second countable), see  Corollary D.1.2 in Appendix D of \cite{CC2}. Further, since $\tilde{\F}_{1}$ and  $\tilde{\F}^{\pitchfork}_{1}$ have global product structure,  it follows that $\FF_{1}$ is a simply connected. So in order to prove this proposition we only need to show $\FF_{1}$ is Hausdorff.

If not, there exist $F_{1}$ and $F_{2}$, which are representatives in two different equivalence classes in $\FF_{1}$, such that any pair of open neighborhoods of $F_{1}$ and $F_{2}$ will intersect. We might as well take two small enough open neighborhoods of them, call them $U_{1}$ and $U_{2}$. Then pick two arcs $J_{1}$ and $J_{2}$ transverse to $F_{1}$ and $F_{2}$ in $x$ and $y$ respectively, such that $\pi'(J_{1})=U_{1}$ and $\pi'(J_{2})=U_{2}$, where $\pi':\tilde{M}\rightarrow \tilde{M}\diagup \tilde{\F}_{1}$   is the natural projection, and $\tilde{\F}^{\pitchfork}_{1}(x)\cap F_{2}=y$. If necessary, shrinking $J_{1}$ and $J_{2}$. Since $U_{1}\cap U_{2}\neq \varnothing$, there is a leaf $F'$ in $\tilde{\F}_{1}$ intersects the $J_{1}$ and $J_{2}$. Hence, $F'$ will intersect $\tilde{\F}^{\pitchfork}_{1}(x)$ twice, because $\F_{1}$ is continuous and $\tilde{\F}_{1}$  transverse to $\tilde{\F}^{\pitchfork}_{1}$.  But this will contradict the property of global product structure by Theorem \ref{gps}.

\end{proof}
\end{prop}

\begin{rmk}\label{rega}
In fact, by the argument and method in the preceding proof,  $\tilde{\F}^{\pitchfork}_{1}(x)$ can be regarded as the leaf space $\FF_{1}$ of $\tilde{\F}_{1}$ for any $x\in \tilde{M}$.
\end{rmk}

\begin{prop}\label{tror}
Suppose $\F_{1}$ is a codimension one foliation without holonomy on the closed manifold $M$. Then $\F_{1}$ is transversely orientable.
  \begin{proof}
We give a concise proof when $\F_{1}$ is a hyperplanar foliation, for the general case see Theorem VIII.2.2.8 of \cite{HH2}. If the hyperplanar foliation $\F_{1}$ is not transversely orientable, then any foliation $\F^{\pitchfork}_{1}$ transverse to $\F_{1}$ will be not orientable. Take a leaf $\tilde{\F}^{\pitchfork}_{1}(x)$, $x\in \tilde{M}$, as the leaf space of $\tilde{\F}_{1}$, then by the global product structure of $\tilde{\F}_{1}$ and $\tilde{\F}^{\pitchfork}_{1}$, we can define the homeomorphisms on $\tilde{\F}^{\pitchfork}_{1}(x)$ induced by the deck transformations $\pi_{1}(M)$ as follows: for any $\alpha\in\pi_{1}(M)$,

\begin{equation*}
\begin{aligned}
     \hat{\alpha}&  : &\tilde{\F}^{\pitchfork}_{1}(x) & \rightarrow & & \tilde{\F}^{\pitchfork}_{1}(x)\\
           & & y           & \mapsto & & \hat{\alpha}(y)=\alpha(\tilde{\F}_{1}(y))\cap\tilde{\F}^{\pitchfork}_{1}(x)
\end{aligned}
\end{equation*}

Since each leaf of $\F_{1}$ is homeomorphic to the hyperplane, then any homeomorphism $\hat{\alpha}$ induced by $\alpha\in\pi_{1}(M)\diagdown\{0\}$ is fixed point free. Since   $\F^{\pitchfork}_{1}$ is not orientable, there exists an element $\beta$ in $\pi_{1}(M)$ such that the induced homeomorphism $\hat{\beta}  : \tilde{\F}^{\pitchfork}_{1}(x)  \rightarrow   \tilde{\F}^{\pitchfork}_{1}(x)$ reverses the orientation, then  $\hat{\beta}$ has a fixed point, which is a contradiction.

  \end{proof}
\end{prop}

\subsection{Hyperplanar foliations}\label{sec2.2}

In this subsection, we will study the hyperplanar foliations. Recall that a codimension one foliation on an $n$-dimensional closed manifold  is said to be a \emph{hyperplanar foliation} if all its leaves are homeomorphic to $\R^{n-1}$, denote it as $\F$. And the hyperplanar foliations considered in this subsection are also assumed to be  $C^{0}$. We remark that there are many such foliations.  Typically,  they appear in the condimension one Anosov systems, see\cite{Sma}.

The study of the limitations and effects of foliations on the topology of the ambient manifolds has been in progress since the last century, see\cite{No}\cite{Thu1}\cite{Law}. However, it is extremely difficult to investigate the general foliations. The known elegant results appear in the case of low-dimensional manifolds, and the foliations considered there are codimension one and have high regularity, like $C^{r}(r\geq 2)$ or analytic. For example, in \cite{Ha} H. Rosenberg proved that the torus $\T^{3}$ is the only closed smooth manifold which admits a $C^{2}$ codimension one  foliation by planes. W. Thurston pointed out in \cite{Thu}, a compact manifold admits a codimension one foliations if and only if the Euler-Poincar\'{e} characteristic is zero. C.  Palmeira showed that a simply connected ($n+1$)-manifold, foliated by $n$-planes (all the leaves are diffeomorphic to the standard $\R^{n}$), $n\geq 2$, is determined up to conjugacy by the leaf space in \cite{Pal}. It implies that,  the only simply connected noncompact manifolds that admit such
foliations are, up to diffeomorphism, the Euclidean spaces $\R^{n+1}$.

As expected, with more restrictive hypothesis, the relevant results will be more clear. When the object under consideration is the hyperplanar foliations, combining the discussion about foliations without holonomy in the previous subsection, we should be able to get more topological properties of the ambient manifold. The following several propositions describe part of them.

\begin{prop}\label{unisp}
If $M^{n}$ admits a hyperplanar foliation $\F$, then $\tilde{M}$ is homeomorphic to $\R^{n}$.
\begin{proof}
By  Theorem \ref{existrans}, we can choose a foliation $\F^{\pitchfork}$ transverse to $\F$. Since all the leaves of $\F$ are homeomorphic to $\R^{n-1}$, their holonomy groups are trivial, then $\F^{\pitchfork}$ and $\F$ have global product structure follows from  Theorem \ref{gps}. It implies that each leaf of $\tilde{\F}^{\pitchfork}$ is homeomorphic to $\R$. Therefore, Proposition \ref{glho} concludes the proof.

\end{proof}

\end{prop}

And using the foliation $\F$ we are also able to show:

\begin{prop}\label{funda}
If $M^{n}$ admits a hyperplanar foliation $\F$, then the fundamental group $\pi_{1}(M)$ is $\Z^{n}$.

\begin{proof}
 $\pi_{1}(M)$ acts as a group of homeomorphisms of $\FF$, the leaf space of $\F$, which is well defined by the property of  global product structure and $\pi_{1}(M)$ preserves $\tilde{\F}$. These nontrivial homeomorphisms are fixed point free, otherwise if there is $\alpha \in \pi_{1}(M)$ and  $\alpha \neq 0$ such that $\alpha(F)=F$ for some leaf $F$ in $\tilde{\F}$, then $\pi(F)=\pi(\alpha(F))$ implies that there exists a non-contractible loop in $\pi(F)$, which is impossible. Hence $\pi_{1}(M)$ acting freely on $\FF\simeq \R$, by  H\"{o}lder theorem (see \cite{CC} or \cite{HH2}),  $\pi_{1}(M)$ is free Abelian.

By the theorem of the structure of finite generated groups, we can get $\pi_{1}(M)\cong \mathbb{Z}^{k}$, $k\geq 0$, then $k$ must be equal to $n$, which follows from the cohomology theory of closed manifolds.

\end{proof}

\end{prop}

\begin{rmk}
In fact, if the foliation $\F$ is induced by the stable distribution of a partially hyperbolic diffeomorphism $f$, $\F$ is invariant under $f$. And there will be another proof of the proposition above.
Indeed, for any leaf $F$ of $\tilde{\F}$, the $\pi:\tilde{M}\rightarrow M$ restricted to $F$ is injective. Since $\pi$ is a local homeomorphism, there exists $\epsilon>0$ such that $0<d_{\tilde{M}}(x,y)<\epsilon$, where $d_{\tilde{M}}$ is lifted from a Riemannian metric on $M$. Assume $x,\ y\in F$ and $\pi(x)=\pi(y)$, then for all $i\geq 0$, $\pi(\tilde{f}^{i}(x))=\pi(\tilde{f}^{i}(y))$. However, for large enough $n$, one can have $d_{\tilde{M}}(\tilde{f}^{n}(x),\tilde{f}^{n}(y))<\epsilon$, which is a contradiction.  So, $\pi_{1}(M)\diagdown\{0\}$ must act without fixed point on $\FF$. If not, there exists a $\alpha\in\pi_{1}(M)\diagdown\{0\}$ such that $\alpha(F_{0})=F_{0}$ for some leaf $F_{0}$ of $\tilde{\F}$, then for any point $x\in F_{0}$, we have $\pi(x)=\pi(\alpha(x))$, which contradicts the injectivity of $\pi\mid_{\tilde{\F}}$.
\end{rmk}

Finally, we are able to prove the following useful proposition, which can be regarded as a generalization of Rosenberg's result \cite{Ha}.

\begin{prop}\label{torus}
  If a smooth closed manifold $M^{n}$ admits a $C^{0}$ hyperplanar foliation, then $M^{n}$ is homeomorphic to $\T^{n}$.
\begin{proof}

When $n=2$, $M^{2}$ is $\mathbb{T}^{2}$ by the theorem of classification of closed surfaces.

When $n=3$, we  can invoke a  powerful theorem about the classification of orientable closed 3-manifolds to avoid the usual tedious analysis. Note that universal covering of $M^{3}$ is $\mathbb{R}^{3}$, then $M^{3}$ is irreducible, so $M^{3}$ is prime (\cite{Hat}, Proposition 1.4). Hence,  $M\cong\mathbb{T}^{3}\# S^{3}\cong\mathbb{T}^{3}$ by Theorem 2.1.3 in \cite{As}.

When $n=4$, the classifying map $\phi:M^{4} \rightarrow \mathbb{T}^{4}$ is a homotopy equivalence, since $\mathbb{T}^{4}$ is a $K(\mathbb{Z}^{4},1)$-space. Then by the main theorem of section 11.5 of Freedman and Quinn \cite{FQ} (this theorem shows the classification of aspherical manifolds with poly-(finite or cyclic) fundamental group), $\phi$ is homotopic to a homeomorphism, because all finitely generated abelian groups are polycyclic and $M^{4}$ is a $K(\mathbb{Z}^{4},1)$-space, so it is aspherical. In other words, $\mathbb{T}^{4}$ is topologically unique in $K(\mathbb{Z}^{4},1)$-spaces.

When $n\geq 5$, since $M^{n}$ is homotopic equivalent to $\mathbb{T}^{n}$. Due to Hsiang and Wall \cite{HW}, any closed $n$-manifold homotopy equivalent to a torus is homeomorphic to a torus.

Therefore, the desired conclusion holds.

\end{proof}
\end{prop}

It is easy to see that item $(i)$ of Theorem \ref{final} follows from the above proposition.

\subsection{The quasi-isometry }

Moreover, there are some more specific topological characterizations about the hyperplanar foliation $\F$. In this subsection, we provide some results in this regard, which are useful for the following proofs. The hyperplanar foliation $\F$ we are considering in this subsection  is assumed to be $C^{0,1}$.

\begin{prop}\label{nullho}
There is no closed loop transverse to $\F$ which is nullhomotopic.
\begin{proof}
 Suppose the contrary. If there exists a closed contractible transversal $l$, by Theorem \ref{existrans}, we can construct a foliation $\F^{\pitchfork}$ transverse $\F$ on $\T^{n}$ with $l$ as its a closed leaf. Since $\F$ is a codimension one foliation without holonomy, it follows that $\F$ and $\F^{\pitchfork}$ have global product structure. However, $l$ is nullhomotopic, so $\tilde{l}$ is also a closed loop and its   lift on $\R^{n}$  cannot intersect the certain leaf of $\tilde{\F}$ only once, which contradicts the definition of global product structure. This concludes the proof.

\end{proof}

\end{prop}

We say a leaf  $\tilde{\F}(y)$ of $\tilde{\F}$ is \emph{properly embedded} if the embedding $i:\tilde{\F}(y)\rightarrow \R^{n}$ is a proper map, $i.e.$, for any $x,\ x_{n}\in \tilde{\F}(y)$ and $d_{\tilde{\F}}(x,x_{n})\rightarrow \infty$ implies that $d(x,x_{n})\rightarrow \infty$. In other words, the intersection of $\tilde{\F}(y)$ with any compact $K\subset \R^{n}$ is compact.

\begin{prop}\label{embed}
All the leaves of $\tilde{\F}$ are properly embedded $C^{1}$-copies of $\R^{n-1}$. Furthermore, there exists $\delta>0$ such that every Euclidean ball of radius $\delta$ can be covered by a coordinate neighborhood $U$ such that the intersection of each leaf $F$  with $U$ is either empty or a connected set being a hyperplanar disk.

\begin{proof}
Assume that there is  a leaf $F$ of $\tilde{\F}$ is not properly embedded  in $\R^{n}$, namely,  $i:F\rightarrow \R^{n}$ is not proper. Then there exists a compact set $U\subseteq \R^{n}$ such that $V=i^{-1}(U)\cap F$ is not compact. Hence we can find a sequence of points $\{x_{j}\}\subseteq V$ such that $\|x_{j}\|\rightarrow \infty\ as\ j\rightarrow\infty$, but $i(\{x_{j}\})$ accumulate on a point $\bar{x}\in U$. Since the foliations $\tilde{\F}$ and $\tilde{\F}^{\pitchfork}$ are transverse, $\tilde{\F}^{\pitchfork}(\bar{x})$ would intersect $F$ an infinite number of times. There will be a contradiction.

As for the latter statement in this proposition, just notice that there is no $C^{1}$ curve transverse to $\tilde{\F}$ that intersects the $F$ more than once by Proposition \ref{nullho}. Then it follows from the same method as in Lemma 3.2 in \cite{BBI}.

\end{proof}

\end{prop}

\begin{rmk}
Note that under the premise that Proposition \ref{embed} holds, the $\mathbb{R}^{n}$ can be separated by $\tilde{\F}(y)$ into two components for any $y\in\R^{n}$. Indeed, we can view the homeomorphism $h:\tilde{\F}(y)\rightarrow\R^{n}$ as a map $h:\R^{n-1}\rightarrow \R^{n}$. Since $h$ is proper, $h$ can be extended to the one-point compactification, so we can regard it as an embedding $h:\bS^{n-1}\rightarrow \bS^{n}$. Hence, by  Theorem 2B.1 in Hatcher's book \cite{Al}, $\bS^{n}\backslash h(\bS^{n-1})$ has two connected components, which implies that $\R^{n}\backslash \tilde{\F}(y)$ has two connected components, and $\tilde{\F}(y)$ is the boundary of the connected components.

\end{rmk}

By the same argument, the leaves of $\tilde{\F}^{\pitchfork}$ are also properly embedded as well, but we can say something further about $\tilde{\F}^{\pitchfork}$, namely, there exists a  ``length versus volume" inequality for neighbourhoods of leaves of $\tilde{\F}^{\pitchfork}$.

\begin{prop}\label{lvv}
There is a constant $C$ such that for every segment $L$ of a leaf of $\tilde{\F}^{\pitchfork}$, the following inequality
\begin{center}
  $\emph{Volume}(U_{1}(L))\geq C\cdot \emph{Length}(L)$
\end{center}
holds, where $U_{1}(L)\triangleq \{x\in \R^{n}| dist(x,\ L)< 1\}$. 
\begin{proof}

 By the compactness of $\T^{n}$, we can take some flow boxes of radius $2\delta$ to cover $\R^{n}$.  And by  Proposition \ref{embed}, we may choose a constant $l>0$, such that for each segment of a leaf of $\tilde{\F}^{\pitchfork}$  with length $l$, its volume of $1$-neighborhood is not less than $\delta$. Decompose $L$ into some disjoint subsegments $J_{1},\ J_{2},\ \ldots,\ J_{n}\subset L$ of length $l$. Then we have
\begin{center}
  $\text{Volume}(U_{1}(L))\geq \lceil \frac{\text{Length}(L)}{l}\rceil \times \delta$.
\end{center}
In fact, $U_{1}(J_{1})$, $U_{1}(J_{2})$, $\ldots$, $U_{1}(J_{n})$ are pairwise disjoint. Otherwise, we can find a subsegment $\gamma$ of $\tilde{\F}^{\pitchfork}$ with the distance between the endpoints of $\gamma$ is less than $\delta$. Then perturbing such a segment yields a transverse nullhomotopic loop, which contradicts  Proposition \ref{nullho}.

\end{proof}

\end{prop}

Compared with the results about classification of foliations by planes in \cite{Pal}, the foliation $\tilde{\F}$ considered here  should have better geometric properties and more rigidity, because it is obtained from the lift of a foliation on $\T^{n}$. Now, we declare a such property of $\tilde{\F}$.

\begin{prop}\label{bound}
Let $\F$ be a hyperplanar foliation of  $\T^{n}$. Then there exists a linear hyperplane 
$P\subset \R^{n}$ and $T>0$ such that, every leaf of $\tilde{\F}$ lies in a $T$-neighborhood of translate of $P$  and for every $x\in \R^{n}$ the $T$-neighborhood of leaf  $\tilde{\F}(x)$ contains $P+x$, a translation of the linear hyperplane $P\subset\R^{n}$. Furthermore, the linear hyperplane $P$ is unique up to translation.

\end{prop}

The case of three dimension of proposition above can be found in \cite{RP1}, and their proofs are inspired by the arguments in the proof of Theorem 1.3 in \cite{BBI}. Fortunately, in the light of the existence of global product structure, we can give a more concise proof.

\noindent\begin{proof3}
Firstly, we give some definitions and notations, and some necessary facts.

\begin{figure}[htbp]
		\centering
		\includegraphics[width=10cm]{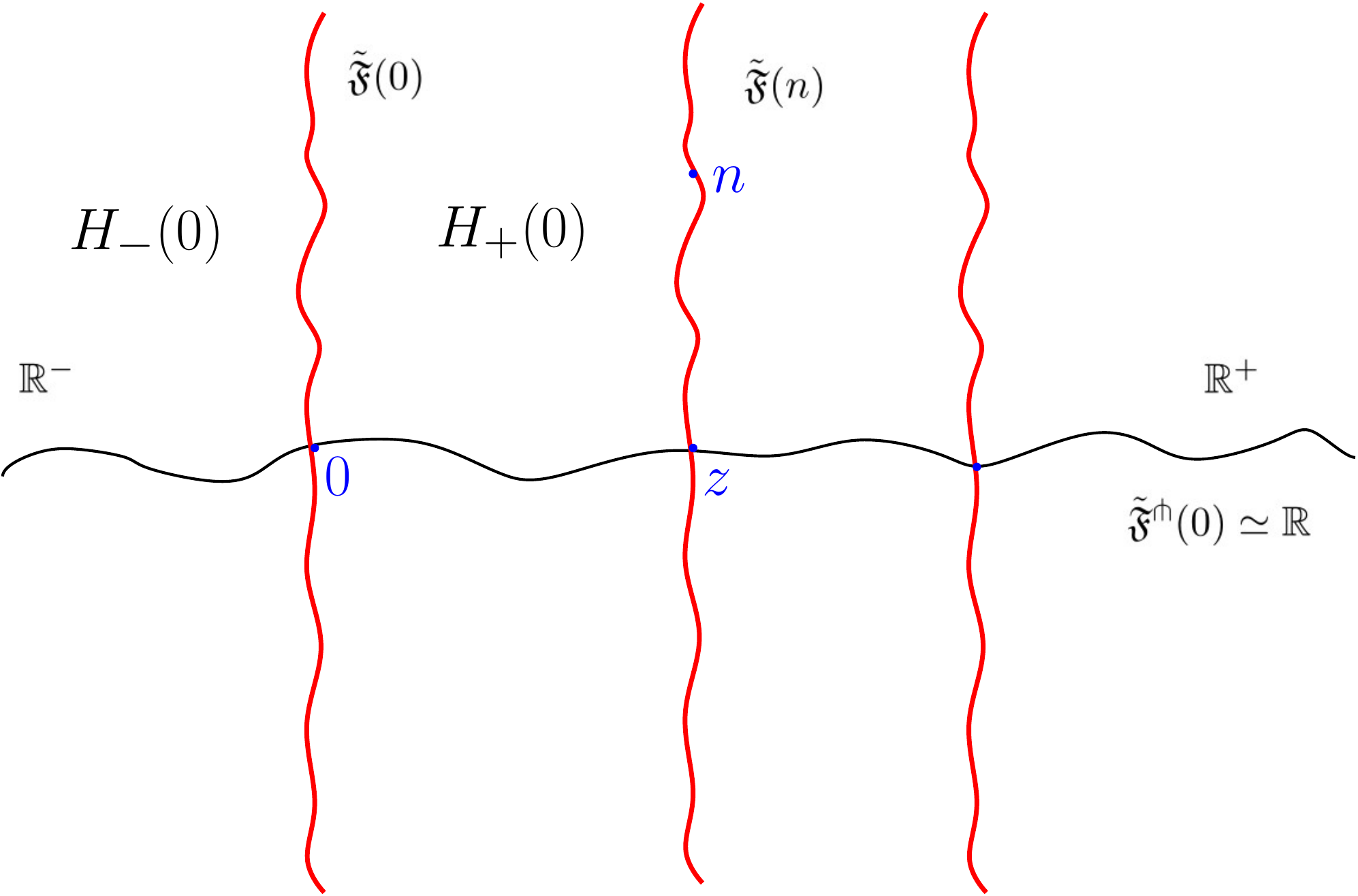}
		\caption{Action on the leaf space of $\tilde{\F}$}
        \label{fig:plane}
\end{figure}

By Remark \ref{rega}, choose the leaf $\tilde{\F}^{\pitchfork}(0)\subseteq\tilde{\F}^{\pitchfork}$ contains $0$ in $\R^{n}$ as the leaf space of $\tilde{\F}$. Since $\tilde{\F}^{\pitchfork}$ and $\tilde{\F}$ have global product structure and $\pi_{1}(\T^{n})$ preserves these foliations, we can define the following map on $\tilde{\F}^{\pitchfork}(0)\simeq\R$, which is induced by $\pi_{1}(\T^{n})=\Z^{n}$,

\begin{equation*}
\begin{aligned}
     \alpha&  : & \Z^{n}\times\R & \rightarrow & & \R\\
           & & (k,t)           & \mapsto & & \alpha(k,t)=\alpha_{k}(t)=\tilde{\F}(t+k)\cap\tilde{\F}^{\pitchfork}(0)
\end{aligned}
\end{equation*}

\noindent And for any $k\in\Z^{n}\setminus\{0\}$, the action $\alpha(k,\cdot)$ is a fixed point free homeomorphism on $\R$ and these actions are free Abelian, which follows from the proof of  Proposition \ref{funda}.

Since $\R^{n}$ is simply connected, we can consider an orientation on $\tilde{\F}^{\pitchfork}$ (because $\F$ is transversely orientable by Proposition \ref{tror}, then $\F^{\pitchfork}$ is orientable, and this orientation is preserved under covering transformations).  $\tilde{\F}^{\pitchfork}(0)\backslash \{0\}$ has two connected components which we call $\R^{-}$ and $\R^{+}$ according to the chosen orientation of $\tilde{\F}^{\pitchfork}$ (Theorem \ref{gps} shows, it is impossible that  there is only one connected component. Of course, it also illustrates that each leaf of $\tilde{\F}^{\pitchfork}$ is homeomorphic to $\R$). By Proposition \ref{embed}, denote two components of $\R^{n}$ separated by $\tilde{\F}(0)$ as $H_{+}(0)$ and $H_{-}(0)$ depending on whether they contain $\R^{+}$ and $\R^{-}$, abbreviated as $H_{+}$ and $H_{-}$ hereinafter, as shown in  Figure \ref{fig:plane} above. We remark that both $H_{+}$ and $H_{-}$ contain the boundary $\tilde{\F}(0)$.

We consider the following subsets of $\Z^{n}$ seen as deck transformations:
\begin{center}
$\Gamma_{+}\triangleq\{k\in\Z^{n}|\alpha_{k}(0)\in H_{+}\},$
\end{center}

\begin{center}
$\Gamma_{-}\triangleq\{k\in\Z^{n}|\alpha_{k}(0)\in H_{-}\}.$
\end{center}

\noindent And because of global product structure we can get
$\Gamma_{+}\cup\Gamma_{-}=\Z^{n}$. On the other hand, since  $\alpha_{k}(\cdot)$ is fixed point free, we have $\Gamma_{+}\cap\Gamma_{-}=\{0\}$.

Further, by the H\"{o}lder theorem, the each action $\alpha_{k}(\cdot)$ is conjugate to an action by translation on $\tilde{\F}^{\pitchfork}(0)\simeq\R$, then for every $k\in\Z^{n}$ we can obtain the following facts:
\begin{center}
$\mbox{for}\ k\in\Gamma_{+},\ \alpha_{k}(\R^{+})\subseteq\R^{+}\ \mbox{and}\  H_{+}+k\subseteq H_{+},$
\end{center}

\begin{center}
$\mbox{for}\ k\in\Gamma_{-},\ \R^{+}\subseteq\alpha_{k}(\R^{+})\ \mbox{and}\ H_{+} \subseteq H_{+}+k.$
\end{center}
\noindent Observe that deck transformations preserve the orientation and the foliation $\tilde{\F}$, then for any $y\in \R^{n}$ and $k\in \Z^{n}$, $H_{\pm}(y)+k=H_{\pm}(y+k)$.

Then, define $\Delta_{+}=\bigcap\limits_{k\in \mathbb{Z}^{n}}\overline{H_{+}+k}$ and $\Delta_{-}$ in a similar way.

\begin{af}\label{nothing}
$\Delta_{+}=\Delta_{-}=\varnothing$.
\end{af}

\noindent\emph{Proof of the Claim 1: } Assume the negation holds, and take $x\in \Delta_{+}$ (we just consider the case $\Delta_{+}\neq\varnothing$ since the  $\Delta_{-}\neq\varnothing$ is symmetric). By the definition of $\Delta_{+}$, $x$ can not be contained in $\tilde{\F}(0)$, so $\tilde{\F}(0)$ and $\tilde{\F}(x)$ are different leaves of $\tilde{\F}$, $\tilde{\F}(x)$ do not intersect with $\tilde{\F}(0)$. Thus $\tilde{\F}(x)$ only lies in the one of $H_{+}$ and $H_{-}$. Without loss of generality, assume $\tilde{\F}(x)\subseteq H_{-}$. The global product structure implies that $\tilde{\F}^{\pitchfork}(0)\cap\tilde{\F}(x)$ has a unique point, denote it $z$. From the proof of the H\"{o}lder theorem, we can know that $\Z^{n}$ induces an Archimedean order $\preceq$ (a bi-invariant order $\preceq$ on a group $G$ is said to be \emph{Archimedean} if for all $g$ and $h$ in $G$ such that $g\neq id$, there exists $n\in\Z$ satisfying $h\prec g^{n}$, more details see \cite{HH2}) on $\tilde{\F}^{\pitchfork}(0)$, then there exists $k'\in\Z^{n}$ such that $k'(\tilde{\F}(0))\cap\tilde{\F}^{\pitchfork}(0)\prec z$. Thus, $\tilde{\F}(x)\in \overline{H_{+}+k'}$, namely, $\tilde{\F}(x)\notin (H_{+}+k')\cap H_{+}$. Since
$\tilde{\F}(x)\in \overline{H_{+}(0)}^{c}$, then $x\notin\Delta_{+}$. Hence the contradiction  concludes this proof of the claim.

\hfill $\blacksquare$\\

In fact, adapting the  same argument in the proof of the claim above, we can prove $\Delta_{+}(y)=\Delta_{-}(y)=\varnothing$ for each $y\in \R^{n}$, where $\Delta_{\pm}(y)=\bigcap\limits_{k\in \mathbb{Z}^{n}}\overline{H_{\pm}(y)+k}$. Therefore, for each point $y\in\R^{n}$, we have that

\begin{equation}
\bigcup\limits_{k \in \mathbb{Z}^{n}}\left(H_{+}(y)+k\right)=\bigcup\limits_{k \in \mathbb{Z}^{n}}\left(H_{-}(y)+k\right)=\R^{n}. 
\nonumber
\end{equation}

Next, we borrow the method  from the proof of Lemma $3.12$ of \cite{BBI} (and the argument after that lemma) to prove the following key claim.

\begin{af}
$\Gamma_{+}$ and $\Gamma_{-}$ are half lattices (this means that there exists a linear hyperplane $P \subset \mathbb{R}^{n}$ such that each one is contained in a half space bounded by $P$). And for every $y \in \mathbb{R}^{n}$ there exists a linear hyperplane $P(y)$ such that $H_{+}(y)$ lies in a half space bounded by $P_{+}(y)$ and $H_{-}(y)$ lies in a half space bounded by $P_{-}(y)$, where $P_{+}(y)$ and $P_{-}(y)$ are linear hyperplanes parallel to $P$.
\end{af}

\noindent\emph{Proof of the Claim 2: }Consider the convex hulls of $\Gamma_{+}$ and $\Gamma_{-}$. If their interiors are disjoint, they can be separated by a linear hyperplane because they are nonempty open convex sets. Then this linear hyperplane is a desired $P$. Otherwise, when the intersection of two convex hulls has nonempty interior, one can consider a point whose coordinates are non-zero rational. Then it is positive rational convex combinations of vectors in $\Gamma_{+}$ as well as of vectors in $\Gamma_{-}$. It follows that $\Gamma_{+} \cap \Gamma_{-}$ has rank not less than 1, which is impossible by the statements below  Proposition \ref{bound}. This implies that there is a linear hyperplane $P$ separating the above convex hulls.

Thus, the first half of the claim is proved. As for the latter part, the proof is as follows.

Take $z \in \R^{n}$ and let $\co_{+}(z)\triangleq\left(z+\Z^{n}\right) \cap H_{+}$. We have that $\mathcal{O}_{+}(z) \neq \varnothing$ (otherwise $z \in \Delta_{-}$ contradicting  Claim \ref{nothing}). Moreover, since $\Gamma_{+}$ preserves $H_{+}$, we also have that $\co_{+}(z)+\Gamma_{+} \subset \co_{+}(z)$. The symmetric statements hold for $\co_{-}(z)\triangleq\left(z+\Z^{n}\right) \cap H_{-}$.

Since $\Gamma_{+}$ and $\Gamma_{-}$ are half-lattices separated by a linear hyperplane $P$, it is easy to know that $\co_{+}(z)$ and $\co_{-}(z)$ are separated by a linear hyperplane $P_{z}$ parallel to $P$. Then choose a $\delta$ given by Proposition \ref{embed} such that every point $z$ has a neighborhood $U_{z}$ containing $B_{\delta}(z)$ and such that $\tilde{\F}(x) \cap U_{z}$ is connected for every $x \in U_{z}$. Let $\left\{z_{i}\right\}_{i=1}^{m}$ be a $\delta / 2$-net in a fundamental domain $\Omega$, namely, $\left\{z_{i}\right\}_{i=1}^{m}$ are $\delta / 2$-dense in  $\Omega$. And let $P_{z_{i}}^{+}$ and $P_{z_{i}}^{-}$ denote the half spaces defined by the linear hyperplanes $P_{z_{i}}$ parallel to $P$ containing $\co_{+}\left(z_{i}\right)$ and $\co_{-}\left(z_{i}\right)$, respectively. Then, $H_{+}$ is contained in the $\delta$-neighborhood of $\bigcup_{i} P_{z_{i}}^{+}$and the symmetric statement holds for $H_{-}$.
This  implies that $\tilde{\F}(0)$ is contained in the intersection of these half spaces, which is a strip bounded by linear hyperplanes $P_{+}$ and $P_{-}$ parallel to $P$.

Hence, we have proved the whole claim.

\hfill $\blacksquare$\\

Note that the basis point $0$ considered in the above proof is not essential. Thus we have proved that, for every $y \in \mathbb{R}^{n}$ there exists a linear hyperplane $P(y)$ and translates $P_{+}(y)$ and $P_{-}(y)$ such that $H_{\pm}(y)$ lies in a half space bounded by $P_{\pm}(y)$.

It remains to show that there exists a uniform constant $T$ such that, for any $z\in\R^{n}$ the $T$-neighborhood of $\tilde{\F}(z)$ contains $P_{+}(z)$. Consider $y\in\R^{n}$, and denote the distance between $P_{+}(y)$ and $P_{-}(y)$ as $T(y)$ (and from the previous proof, we know that $\tilde{\F}(y)$ lies at distance smaller than $T(y)$ from $P_{+}(y)$).  In order to obtain it, it is sufficient to prove that the projection from $\tilde{\F}(y)$ to $P_{+}(y)$ by an orthogonal vector to $P(y)$ is surjective. This is obvious, otherwise, there is a segment $l$ joining $P_{+}(y)$ to $P_{-}(y)$, but $l$ does not intersect $\tilde{\F}(y)$, which  is impossible, because  $\tilde{\F}(y)$ lies in the strip bounded by $P_{+}(y)$ and $P_{-}(y)$ and  each curve from $H_{-}(y)$ to $H_{+}(y)$ must intersect $\tilde{\F}(y)$.

Since the foliation $\tilde{\F}$ is invariant under the translations of $\Z^{n}$, then by the compactness of $\T^{n}$,  we obtain that $T(y)$ above can be chosen uniformly bounded. Further, because all the leaves of $\tilde{\F}$ are disjoint from each other, the linear hyperplane $P(y)$ is  parallel to $P$, so $P(y)$ cannot depend on $y$. Therefore, the desired $P$ is unique up to translation.

This concludes the proof of the Proposition \ref{bound}.

\end{proof3}

Next, we will study the hyperplanar  foliation $\F$ on a Riemannian manifold, which will further reveal its shape. More directly, we will show that $\F$ is quasi-isometric (Proposition \ref{iso1}) and quasi-geodesic (Proposition \ref{quageo}).

\begin{defi}\label{quasi}
Let $M$ be a Riemannian manifold, a foliation $\mathcal{K}$ on  $\tilde{M}$ is called \emph{quasi-isometric}, if there are constants $a,  b\in\R$ such that for any two points $x, y\in \tilde{M}$ which lie on the same leaf of $\mathcal{K}$, we have $d_{\mathcal{K}}(x,y)\leq a\cdot d_{\tilde{M}}(x,y)+b$, where $d_{\mathcal{K}}$ is the path distance measured along the leaf induced by the metric $d_{\tilde{M}}$ on $\tilde{M}$.
\end{defi}

Observe that a leaf is properly embedded in $\tilde{M}$ but it does not have to be quasi-isometric. It is easy to enumerate many such examples, for example, the foliation on $\R^2$ lifted from a foliation of $\T^{2}$ with a Reeb component. In effect, the latter is much stronger, because it implies that there is a uniform control over the size of these two distances, even linear control.

In order to prove the quasi-isometric property of $\tilde{\F}$,  we need to reinterpret and analyze these properness of the embedding, in terms of a comparison between the intrinsic metric $d_{\tilde{\F}}$ and the extrinsic metric $d$ in $\R^{n}$. In the next three propositions, we will show that, in fact, this embedding has better properties in our setting.

\begin{prop}\label{uqi}
$\tilde{\F}$ is uniformly properly embedded in $\R^{n}$ in the following sense: there is a monotonic non-decreasing function $\varphi:\R^{+}\rightarrow\R^{+}$ such that for any two points $x,\ y\in\R^{n}$ and $x\thicksim y$ ($i.e.$, $x$ and $y$ are in the same leaf of $\tilde{\F}$), then $d_{\tilde{\F}}(x,y)\leq \varphi(d(x,y))$. In other words, for every $A>0$ there is $B>0$ so that $\forall x,\ y\in \R^{n},\ x\thicksim y$ and $d(x,y)\leq A$, then $d_{\tilde{\F}}(x,y)\leq B$.
\begin{proof}
Suppose the negation holds, then there exists a sequence $x_{i}\thicksim y_{i}\ (w.r.t)\ F_{i}$ and some constant $t>0$ such that $d(x_{i},y_{i})\leq t$ but $d_{\tilde{\F}}(x_{i},y_{i})\rightarrow \infty$, where $F_{i}$ are leaves of $\tilde{\F}$. Since $\T^{n}$ is compact, we can find some elements $\gamma_{i}\in\Z^{n}$ such that $\gamma_{i}(x_{i})(\triangleq x_{i}+\gamma_{i})$ converges up to subsequence, so without loss of generality assume it does converge. Namely, $\gamma_{i}(x_{i})\rightarrow x\in\R^{n}$. Because of $d(x_{i},y_{i})\leq t$, then  $d(\gamma_{i}(x_{i}),\gamma_{i}(y_{i}))\leq t$, which implies that $\gamma_{i}(y_{i})\rightarrow y\in\R^{n}$ up to another subsequence. Thus, $d(x,y)\leq t$.

Suppose $x\thicksim y$, $i.e.$, there is a leaf $F$ of $\tilde{\F}$ such that $x,\ y\in F$, then $d_{\tilde{\F}}(x,y)= m$. Thus $d_{\tilde{\F}}(x_{i},y_{i})\rightarrow m$ as $i\rightarrow \infty$, which is impossible by hypothesis. Indeed, since all the leaves of $\tilde{\F}$ are simply connected, the local product structure means that there is tubular neighborhood $N(l)$ of $l$ which is foliated as product by disks inside the contiguous leaves of $F$. Then it follows that $B^{F_{i}}_{r}(a_{i})\rightarrow B^{F}_{r}(c)$ when $a_{i}\rightarrow c$, where $B^{F}_{r}(a)\triangleq\{ b\in F(a)|d_{\tilde{\F}}(a,b)=r\}$.

Then $x\nsim y$, which implies that there exists a convergent leaf-sequence whose limit is not unique. Thus the leaf space $\F$  of $\tilde{\F}$ is not Hausdorff, but this contradicts that the leaf space of $\tilde{\F}$ is $\R$, see Proposition \ref{leafspace}.

\end{proof}

\end{prop}

\begin{rmk}

In fact, if $\mathfrak{L}$ is a codimension one foliation on a closed manifold $M$, then $\tilde{\mathfrak{L}}$ is uniformly properly embedded on $\tilde{M}$ if and only if the leaf space $\mathcal{L}$ of $\tilde{\mathfrak{L}}$ is Hausdorff. The proof of the sufficiency is the same as the proof of the above proposition. As for the necessity, assume $\mathcal{L}$ is not Hausdorff. Then there are two different leaves $L_{x}$ and $L_{y}$ of $\tilde{\mathfrak{L}}$, and a leaf sequence $\{L_{i}\}\in \mathfrak{L}$ such that $L_{i}$ converges to both $L_{x}$ and $L_{y}$ on compact subsets respectively. Take two sequence ${x_{i}},\ {y_{i}}$ in $\tilde{M}$ and $x_{i}\thicksim y_{i}\ (w.r.t.)\ L_{i}$ such that $x_{i}\rightarrow x\in L_{x}$ and $y_{i}\rightarrow y\in L_{y}$. If every $x_{i}$ can be joined to $y_{i}$ by an arc $l_{i}\subset L_{i}$ of length bounded by a constant  $r$, then $\{l_{i}\}$ converges (up to subsequence) to a rectifiable arc $l$ of length bounded by $r$ joining $x$ to $y$, and it is contained in a single leaf of $\tilde{\mathfrak{L}}$ , which contradicts the choice of $x,\ y$.

\end{rmk}

The property of uniformly properly embedded does not imply quasi-isometric, because the preceding proper function $\varphi$ may not be linear. Fortunately, in the light of  Proposition $\ref{bound}$ we can prove $\tilde{\F}$ is  indeed quasi-isometric as follows.

\begin{prop}\label{iso1}
$\tilde{\F}$ is quasi-isometric.
\begin{proof}
\begin{figure}[htbp]
		\centering
		\includegraphics[width=10.6cm]{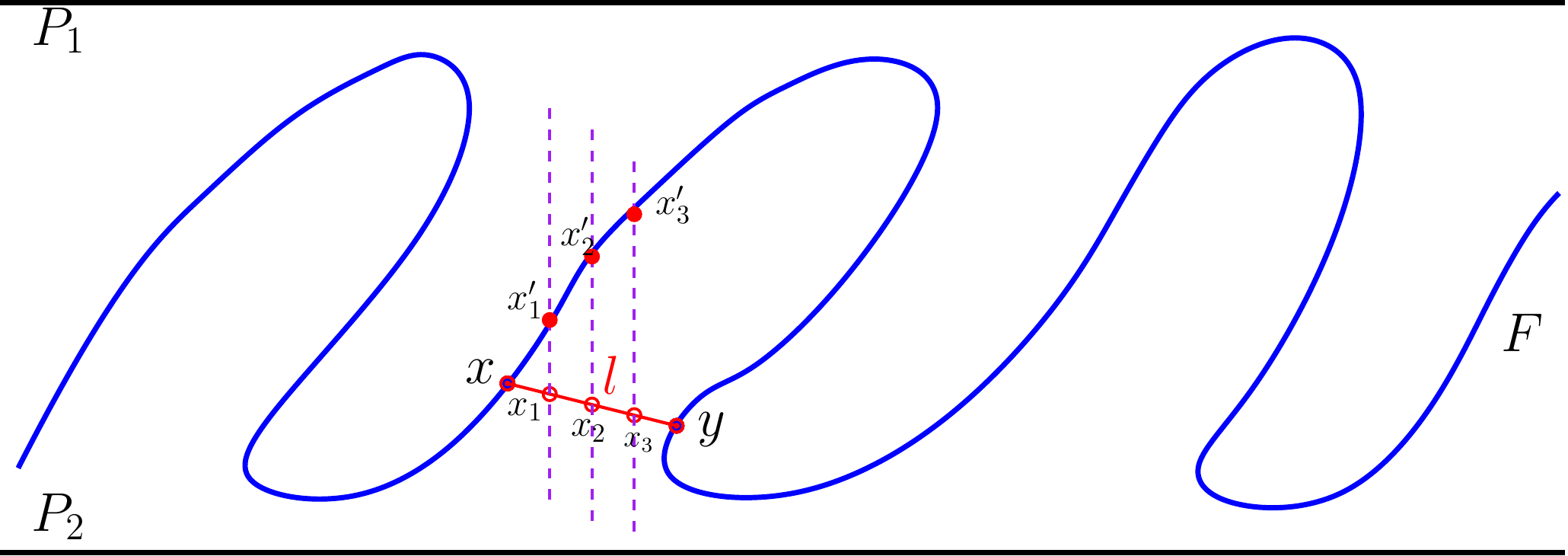}
		\caption{quasi-isometric}
        \label{fig:quaiso}
	\end{figure}
For any leaf $F$ of $\tilde{\F}$, pick any two points $x,\ y\in F$. By Proposition $\ref{bound}$, $F$  lies between two parallel linear hyperplanes $P_{1}$ and $P_{2}$, where $P_{1}$ and $P_{2}$ are translations of the linear hyperplane $P$ and $dist(P_{1},P_{2})=2T$, the $P$ and $T$ are determined by  Proposition $\ref{bound}$. Thus the geodesic $l$ related to $d$ joining $x$ to $y$ also lies between $P_{1}$ and $P_{2}$. Parameterize $l$ by arc length, so that $x=l(0)$ and $y=l(|l|)$. As Figure \ref{fig:quaiso} shows, define $x_{i}\in l$, $0\leq i \leq n$, by $x_{0}=x$, $x_{n}=y$, $d(x_{n-1},x_{n})\leq 1$, $d(x_{i},x_{i+1})=1$, $0\leq i\leq n-2$, where $n \leq |l|< n+1$ (In fact, here we only need to consider $|l|$ can be large enough. If there exist $\{\bar{x}_{i}\}$ and $\{\bar{y}_{i}\}$ in $F$ satisfy that  $d(\bar{x}_{i},\bar{y}_{i})\leq const.$, but $d_{\tilde{\F}}(\bar{x}_{i},\bar{y}_{i})\rightarrow\infty$, then this contradicts the uniformly proper embedding of $\tilde{\F}$, see  Proposition \ref{uqi}). Then for all $0\leq i\leq n$, there is $x_{i}'\in F$ ($x_{0}'=x,\ x_{n}'=y$) such that
\begin{center}
  $d(x_{i},x_{i}')\leq 2T$.
\end{center}
It follows that for all $i\in\{0,1,\cdots,n-1\}$,
\begin{center}
$d(x'_{i}, x'_{i+1})\leq 4T+1$.
\end{center}
By  Proposition \ref{uqi}, $\tilde{\F}$ is uniformly properly embedded with respect to some proper function $\varphi$. Then for all $i$, we have
\begin{center}
  $d_{\tilde{\F}}(x_{i}',x_{i+1}')\leq \varphi(4T+1)$.
\end{center}
Hence,
\begin{center}
  $d_{\tilde{\F}}(x,y)\leq(|l|+1)\varphi(4T+1)=\varphi(4T+1)d(x,y)+\varphi(4T+1)$,
\end{center}
which completes the proof.

\end{proof}

\end{prop}

Furthermore, we may show $\tilde{\F}$ has another interesting geometric feature. We say a foliation $\mathfrak{P}$ on $\tilde{M}$ is \emph{quasi-geodesic} if all of its leaves are uniformly quasi-geodesic, where quasi-geodesic means that for any $x,\ y\in\tilde{M}$, $x\in \mathfrak{P}(y)$, and a minimal length geodesic $l$ joining $x$ to $y$ in $\tilde{M}$, then there is a path $\alpha$ in $\mathfrak{P}(y)$ connecting $x$ and $y$, which is at a bounded distance from $l$. In general quasi-isometric and quasi-geodesic are unrelated, while the foliation $\mathfrak{P}$ is codimension one, then $\mathfrak{P}$ is quasi-geodesic implies that it is also quasi-isometric, see  Theorem 3.2 in \cite{Fen}.

\begin{prop}\label{quageo}
$\tilde{\F}$ is quasi-geodesic.
\begin{proof}

\begin{figure}[htbp]
		\centering
		\includegraphics[width=9cm]{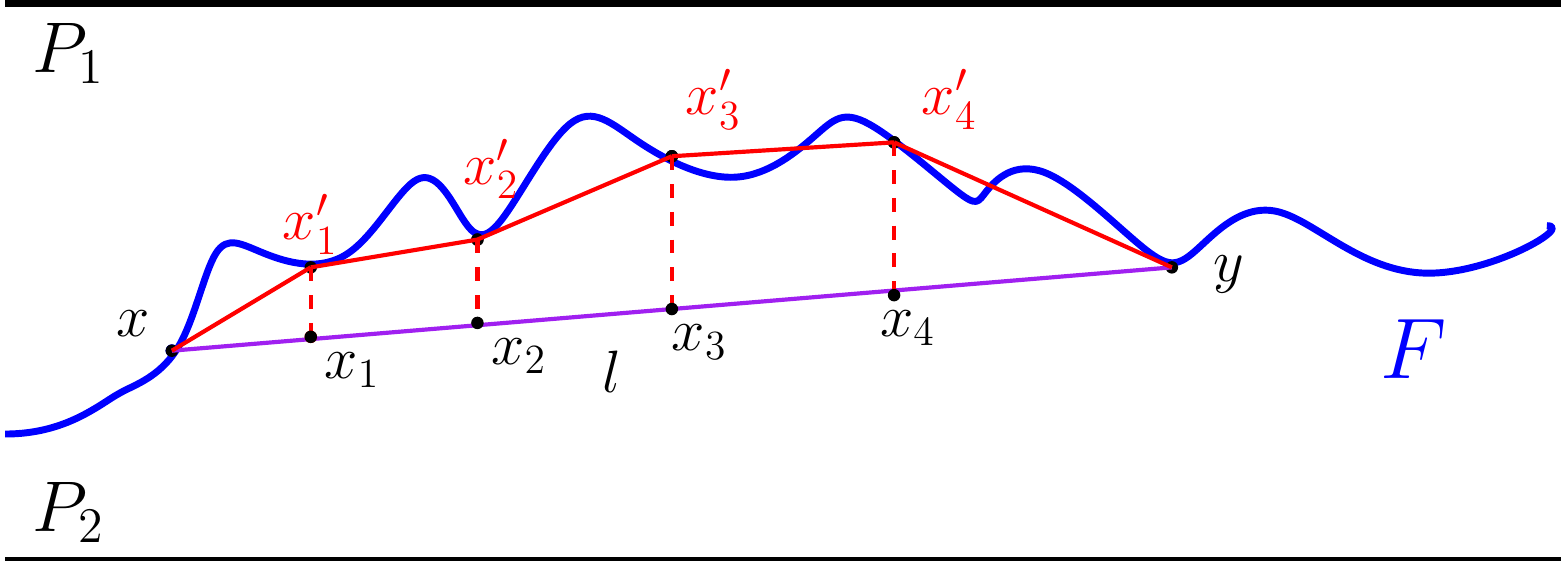}
		\caption{quasigeodesic}
        \label{fig:quasigeo}
\end{figure}

Take arbitrarily a leaf $F$ of $\tilde{\F}^{u}$, and pick any two points $x$ and $y$ in $F$. Then denote $l$ the geodesic joint $x$ and $y$ in $\R^{n}$, and parameterize $l$ by arc length. Choose $x_{i}\in l$ and $x_{i}'\in F$ in the same way as the previous proposition, $0\leq i\leq \lceil |l| \rceil$, as shown in Figure \ref{fig:quasigeo}, where $P_{1}$ and $P_{2}$ are translations of $P$ obtained by Proposition \ref{bound}.

Since $d_{\tilde{\F}}(x_{i}',x_{i+1}')\leq \varphi(4T+1)$ for all $i$, there exists a path $\alpha_{i}\subset F$ connecting $x_{i}'$  and $x_{i+1}'$ contained in a uniformly bounded neighborhood of $l_{i}$ joining $x_{i}'$  and $x_{i+1}'$ in $\R^{n}$. By the selection of $x_{i}'$, we know that all the $l_{i}$ lies in a $2T$-neighborhood of $l$. It follows that there is a path connecting $x$ and $y$ with a bounded distance from $l$. Since $\tilde{\F}$ is lifted from $\T^{n}$, it is invariant under integer translations. Then the compactness implies that distance can be chosen uniformly bounded. Namely, $\tilde{\F}$ is quasi-geodesic.

\end{proof}

\end{prop}

As the end of this section, we show by the way that $\tilde{\F}^{\pitchfork}$ is also quasi-isometric.

\begin{prop}\label{ciso}
$\tilde{\F}^{\pitchfork}$ is quasi-isometric.
\begin{proof}
Let $P$ be the linear hyperplane determined by Proposition \ref{bound}, and pick a unit vector $\upsilon$ orthogonal to $P$ in $\R^{n}$. Since $\tilde{\F}$ and $\tilde{\F}^{\pitchfork}$ have  global product structure, it follows that for every $K>0$ there exists $C$ such that each segment of $\tilde{\F}^{\pitchfork}$ of length $C$ starting at $x$ intersects $P+x+K\upsilon$. Otherwise, there exists $K_{0}>0$ such that for any $n >0$, there is a segment $l_{n}$ of $\tilde{\F}^{\pitchfork}$ of length $n$ starting at $x_{n}$ will not intersect $P+x_{n}+K_{0}\upsilon$. Then by taking a subsequence of $\{l_{n}\}$ and certain translations, we can translate these initial points in a bound region and get a leaf of $\tilde{\F}^{\pitchfork}$ which does not intersect each leaf of $\tilde{\F}$.

The above assertion implies that $\tilde{\F}^{\pitchfork}$ is quasi-isometric. In fact, any segment from $\tilde{\F}^{\pitchfork}$ has a length greater than $tC$.  Then its endpoints are at distance at least $tK$. Hence,
\begin{center}
$\displaystyle
\lceil\frac{d_{\tilde{\F}^{\pitchfork}}(x,y)}{C}\rceil\cdot K\leq d(x,y),
$
\end{center}
which concludes the proof.
\end{proof}
\end{prop}

\bigskip
\section{Locally unique integrability}\label{sec3}

Through the discussion and analysis of the properties of hyperplanar foliations in the previous section, now  we are ready to prove item $(ii)$ of Theorem \ref{final}, that is, the distribution $E^{c}$ is locally uniquely integrable.

First, we show the idea of the proof about the ``absolute" version as described in Remark \ref{rmk0}. Namely, we have

\begin{prop}\label{absolute}
Let $g:N\rightarrow N$ be a codimension one absolutely partially hyperbolic diffeomorphism. Then the distribution $E^{c}$ is locally uniquely integrable. 
\end{prop}

Since $g$ is codimension one, $E^{s}$ is equal to the zero bundle. It is uniquely integrable and the tangent foliation, where each leaf consists of a single point, is trivially quasi-isometric. Thus we can obtain the following lemma from the known result about integrability established by M. Brin, see Theorem 1 of \cite{Bri}.

\begin{lem}\label{cohe}
If $g$ is a codimension one absolutely partially hyperbolic diffeomorphism of a compact n-dimensional Riemannian manifold $M$. Suppose the unstable foliation $\F^{u}$ is quasi-isometric when lifted to the universal covering space, then the distribution $E^{c}$ is locally uniquely integrable.
\end{lem}

Therefore, we can give the following concise proof:

\noindent\begin{proof1}
As a consequence of Proposition \ref{torus}, $N$ must be an $n$-torus $\T^{n}$. Then by Proposition \ref{iso1}, $\F^{u}$ is quasi-isometric. Henceforth, the conclusion follows from Lemma \ref{cohe}.

\end{proof1}

Next,  we give another more direct proof of item $(ii)$ of Theorem \ref{final}, which mainly adopts the methods of \cite{PS}.

\noindent\begin{proof2} As we discussed earlier, $\F^{u}$ is invariant unstable codimension one foliation by hyperplanes on
$N\backsimeq \T^{n}$.  The leaves of $\tilde{\F}^{u}$ are uniformly closed to leaves of a foliation by linear hyperplanes follows from Proposition \ref{bound}. We denote the latter by $\mathcal{P}$, the leaves are linear hyperplanes translated by $P_{0}$, the hyperplane determined by Proposition \ref{bound} and we may assume that it contains the origin of $\R^{n}$.

Replace $f$ with $f^{-1}$ if necessary, we can follow the idea of \cite{PS}. It is worthwhile to indicate that, we will also consider the following case, where the foliation integrated by $(n-1)$-dimensional distribution $E^{s}$ is not complete, as shown in Figure \ref{fig:coexample}.

\begin{figure}[htbp]
		\centering
		\includegraphics[width=7cm]{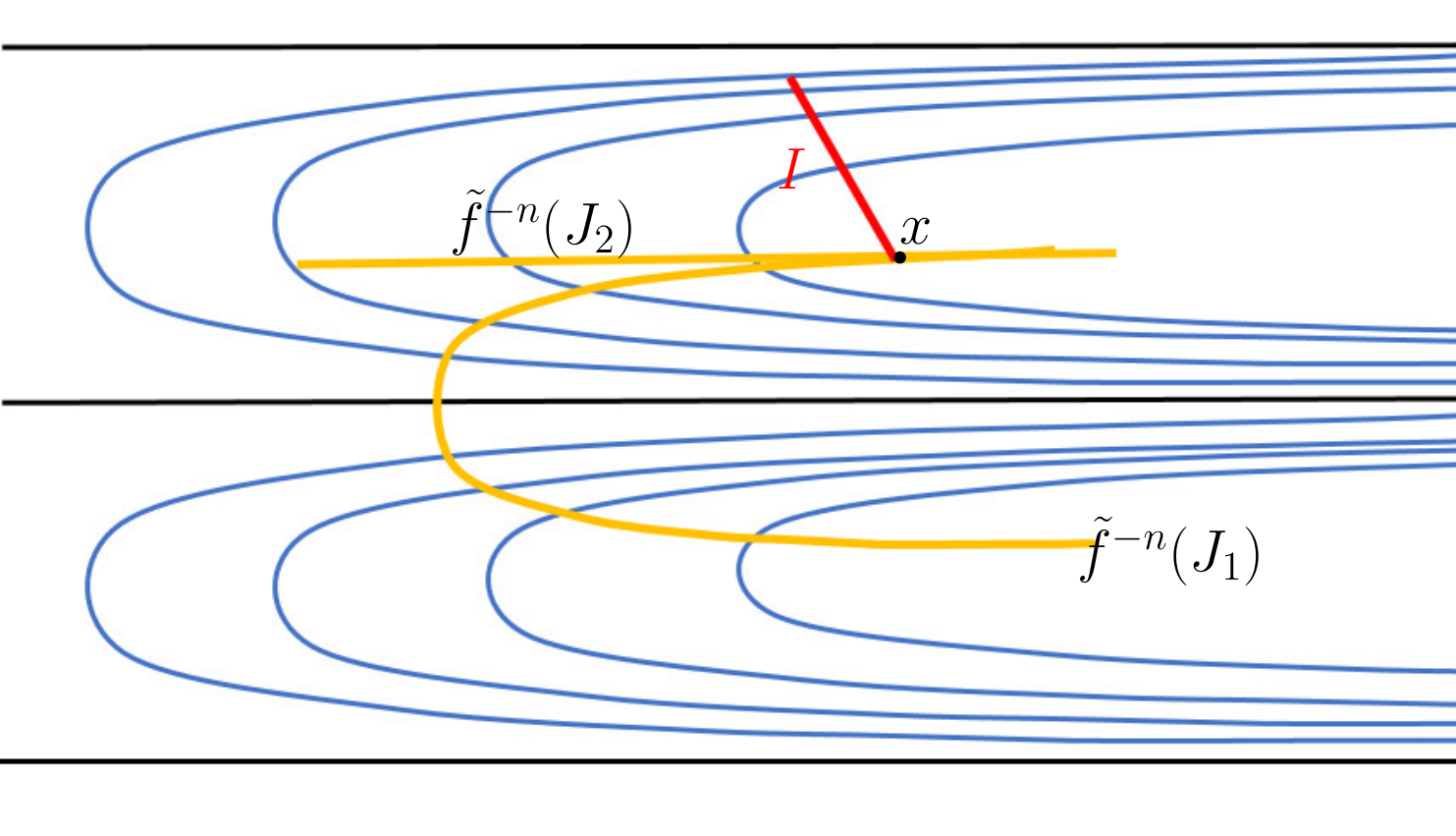}
		\caption{The case involving Reeb components}
        \label{fig:coexample}
\end{figure}

In fact, due to Proposition \ref{bound}, we can exclude the above case. Since for every $x\in \R^{n}$ the $T$-neighborhood of leaf  $\tilde{\F}^{s}(x)$ contains a linear hyperplane $P+x$, then $\tilde{\F}^{s}$ must be complete. Thus, the arguments used in \cite{PS} work for our setting. Finally, we can obtain the  locally unique integrability of $E^{c}$.

For the completeness of the proof, we give the following details of above assertion. To be specific, since the distribution $E^{c}$ is one dimensional, for each $x\in\R^{n}$ there are $C^{1}$ curves passing through $x$ tangent to $E^{c}(x)$ by Peano's Theorem. However, it is not necessarily locally uniquely integrable. As is well known, the regularity of the central distribution is not smooth, even Lipchitz, although $f$ has a better regularity. Next, we will prove that $E^{c}$ in above setting is actually locally uniquely integrable.

Assume that $E^{c}$ is not locally uniquely integrable at some point $x\in \R^{n}$, let $J_{1}$ and $J_{2}$ be two different central segments through $x$, where we call a segment $J$ from integral curves related to $E^{c}$ \emph{central segment}, $i.e.$, for each $x\in J$, then $T_{x}J=E^{c}(x)$, the \emph{stable segments} can be defined in a similar way. We may also assume that $x$ is at the boundary in $J_{1}$ of this intersection, and $|J_{1}|,\ |J_{2}|\leq \delta$, where $\delta$ is determined by the following lemma.

\begin{lem}[\cite{PS}, Lemma 4.3]\label{ps43}
Let $f$ be a  partially hyperbolic $C^{1}$-diffeomorphism with splitting $TN=E^{s}\oplus E^{c}$, where distribution $E^{s}$ is $(n-1)$-dimensional. Then, there is a constant $\delta$ such that for any $\epsilon\leq\delta$ there exists $C=C(\epsilon)$ such that  every central segment $J$ with $|J|\leq \epsilon$ satisfies  $|\tilde{f}^{-n}(J)|\leq C$ for any $n\geq 0$.

\end{lem}

Take $y\in J_{1}\setminus J_{2}$ such that $\tilde{\F}^{s}_{\eta}(y)\cap J_{2}={z}$, where $\tilde{\F}^{s}_{\eta}(y) \triangleq \{ p\in \tilde{\F}^{s}(y)| d_{\tilde{\F}^{s}}(y, p)\geq \eta\}$, and $\eta> 0$ is the radius of balls in $\R^{n}$ such that in every such balls there is no $sc-$bigon ($i.\ e.$, a loop consisting of a stable segment and a central segment). The existence of constant $\eta$ is ensured by the transversality  of $E^{s}$ and $E^{c}$ and the compactness of $N$.

For any $K>0$, there exists $n_{0}$ such that for any $n\geq n_{0}$, we can consider the $K-$cylinder $\tilde{\F}^{s}_{K}(\tilde{f}^{-n}(J_{1}))\triangleq\bigcup\limits_{p\in \tilde{f}^{-n}(J_{1})}\tilde{\F}^{s}_{K}(p)$, which is well defined, since  $\tilde{\F}^{s}_{K}(p)\cap \tilde{\F}^{s}_{K}(q)=\varnothing$ for any $p,\ q\in \tilde{f}^{-n}(J_{1})$ (Otherwise, there is a $sc-$bigon under $\tilde{f}^{n}$  become a smaller $sc-$bigon lies in $B_{\eta}(x)$, which  in contradiction with the choice of $\eta$).  The  set
\begin{center}
$\bigcup\limits_{p\in \tilde{f}^{-n}(J_{1})}\{q\in\tilde{\F}^{s}(p)| d_{\tilde{\F}^{s}}(q,p)=K\}\subset \tilde{\F}^{s}_{K}(\tilde{f}^{-n}(J_{1}))$
\end{center}
is called \emph{$s-$boundary} of that cylinder. For any $L$, there is $K=K(L)$  such that the length of any arc joining $\tilde{f}^{-n}(J_{1})$ with $s-$boundary must be greater than $L$, which because Proposition \ref{bound} implies that the foliation $\tilde{\F}^{s}$ is uniformly bounded by a certain linear foliation.

Let $r= d(y,\ z)$ and $C=C(\delta)$, choose $L\gg C$ and set $K=K(L)$. Then take $n$ large enough such that $q\in \tilde{\F}^{s}_{K}(p)$ implies that $d(\tilde{f}^{n(p)},\tilde{f}^{n}(q))<\frac{r}{2}$. By Lemma \ref{ps43} we can get $|\tilde{f}^{-n}(J_{2})|\leq C$ and $\tilde{f}^{-n}(J_{2})$ does not intersects the $s-$boundary of $K-$cylinder. Hence, $\tilde{f}^{-n}(J_{2})\subset \tilde{\F}^{s}_{K}(\tilde{f}^{-n}(J_{1}))$ and so $\tilde{f}^{-n}(z)\in \tilde{\F}^{s}_{K}(\tilde{f}^{-n}(y))$. Finally we obtain $d(y,z)<\frac{r}{2}$, which contradicts the definition of $r$.

\end{proof2}

\section{The linear part }\label{sec4}

In this section, we will prove the item $(iii)$ of Theorem \ref{final}. Namely, $f$ is semiconjugated to a codimension one Anosov diffeomorphism $f_{\ast}$, where $f_{\ast}$ is the induced linear  transformation of $f$ on $\pi_{1}(\T^{n})$. And for convenience, we still use $f_{\ast}$ to denote its lift on $\R^{n}$.

From now on, let $f:\T^{n}\rightarrow \T^{n}$ be a codimension one  partially hyperbolic diffeomorphism,  and let the linear hyperplane $P_{0}$ be the translation of $P$ in Proposition \ref{bound} and contains the origin $0$. And we denote the foliation generated by  the translation of $P_{0}$ as $\mathcal{P}$.

\begin{prop}\label{inv}
The linear hyperplane $P_{0}$ is invariant under $f_{\ast}$.
\begin{proof}
Notice that $|\tilde{f}(x)-f_{\ast}(x)|\leq const.$ for all $x\in\R^{n}$. Then by Proposition \ref{bound} $\tilde{f}(\tilde{\F}(0))$ lies within a bounded distance from $f_{\ast}(P_{0})$. Therefore, if $f_{\ast}(P_{0})\neq P_{0}$, the leaf $\tilde{f}(\tilde{\F}(0))$ can not be contained in a bounded neighborhood of $P_{0}$, which contradicts with the property of linear hyperplane $P_{0}$ obtained in  Proposition \ref{bound}.

\end{proof}

\end{prop}

The above proposition implies that the linear hyperplane $P_{0}$ is the eigen-subspace generated by the eigenvectors of $f_{\ast}$. We will also use different arguments to gradually obtain more information about the eigenvalues of $f_{\ast}$ in the following three propositions,  their conclusions become stronger in order.

\begin{prop}\label{esti}
The $f_{\ast}$ has an eigenvalue with modulus greater than 1, 
and has another one less than 1.
\begin{proof}
Assume by contradiction that the absolute values of all the eigenvalues of $f_{\ast}$ are less than or equal to 1. Since the volume of balls in $\R^{n}$ has polynomial growth, then the  length of the images of any vector under the iterations of $f_{\ast}$ grows sub-exponentially, and hence so does the diameter of the images of any compact set under the iterations of $\tilde{f}$. Then the images are contained in a sequence of balls whose volume grow sub-exponentially.

We now apply this observation to a disk $D$ in an unstable leaf $\tilde{\F}^{u}(x)$. The $diam_{\tilde{\F}^{u}}(D)$ grows exponentially, but  the images are contained in a sequence of balls with sub-exponential volume growth. This will contradict with the estimation of ``length versus volume", see Proposition \ref{lvv}. Hence, the conclusion holds because the Jacobian determinant of $f_{\ast}$ is 1.

\end{proof}
\end{prop}

\begin{prop}
 The invariant linear hyperplane $P_{0}$ is generated by the eigenvectors corresponding to all eigenvalues whose modulus is not less than 1.
\begin{proof}
Replace $f$ with $f^{-1}$, then the leaves of stable foliation $\tilde{\F}^{s}$ are uniformly close to leaves of $\mathcal{P}$. By Proposition \ref{inv}, we assume by contradiction that there are eigenvectors with the modulus of corresponding eigenvalues greater than 1 in the eigenvectors generating $P_{0}$, they will also produce another invariant linear subspace contained in $P_{0}$, call it $P_{2}$. Then, we can decompose $P_{0}$ as $P_{1}\oplus P_{2}$, similarly decompose $\R^{n}=P_{0}^{\pitchfork}\oplus P_{1}\oplus P_{2}$, where $P_{0}^{\pitchfork}$ is the eigenline transverse to $P_{0}$.

Let $\Pi_{2}:\R^{n}\rightarrow P_{2}$ be the projections defined by the decomposition above, and denote the eigenvector with the smallest modulus in the generated eigenvector of $P_{2}$ as $\beta$, $|\beta| > 1$. Then there is a suitable norm $\parallel\cdot\parallel$ in $\R^{n}$ such that

\begin{center}
$\parallel \Pi_{2}(f_{\ast}(x)-f_{\ast}(y))\parallel \geq |\beta|\cdot\parallel \Pi_{2}(x-y)\parallel$,
\end{center}
for all $x,\ y\in\R^{n}$. By $||\tilde{f}(x)-f_{\ast}(x)||\leq const.$, the above inequality implies that

\begin{center}
$\parallel \Pi_{2}(\tilde{f}^{n}(x)-\tilde{f}^{n}(y))\parallel \geq |\beta|^{n}\cdot(\parallel \Pi_{2}(x-y)\parallel-C)$,
\end{center}
for all $n\geq 1$ and some constant $C$. More specifically, for all $x,\ y\in\R^{n}$ and all $n\geq 1$,\\
\begin{center}
$
\begin{aligned}
& \parallel \Pi_{2}(\tilde{f}^{n}(x)-\tilde{f}^{n}(y))\parallel \\
&\geq \parallel \Pi_{2}(f_{\ast}(\tilde{f}^{n-1})(x)-f_{\ast}
(\tilde{f}^{n-1})(y))\parallel -\\
&\quad \parallel \Pi_{2}(\tilde{f}^{n}(x)-\tilde{f}^{n}(y))-\Pi_{2}
(f_{\ast}(\tilde{f}^{n-1})(x)-f_{\ast}(\tilde{f}^{n-1})(y))\parallel\\
&\geq |\beta|\cdot\parallel \Pi_{2}(\tilde{f}^{n-1}(x)-\tilde{f}^{n-1}(y))\parallel-c\\
&\ldots\\
&\geq  |\beta|^{n}\cdot\parallel \Pi_{2}(x-y)\parallel-c-c\cdot |\beta|- \ldots -c|\beta|^{n-1}\\
&\geq |\beta|^{n}\cdot(\parallel \Pi_{2}(x-y)\parallel-C).
\end{aligned}
$
\end{center}

Since there exists $y_{0}\in\tilde{\F}^{s}(x)$ such that $\parallel \Pi_{2}(x-y_{0})\parallel > C$, then
\begin{center}
$
\parallel\tilde{f}^{n}(x)-\tilde{f}^{n}(y_{0})\parallel\geq const \cdot |\beta|^{n}.
$
\end{center}
Thus there is a contradiction with the fact that every pair of points on the same leaf of the stable foliation $\tilde{\F}^{s}$ will be exponentially close under iterations of $\tilde{f}$.

\end{proof}
\end{prop}

Finally, we can get:

\begin{prop}
The $f_{\ast}$ is Anosov, $i.e.$, the linear transformation $f_{\ast}$ on $\R^{n}$ has no eigenvalue of modulus equals to one.
\begin{proof}
Adapt the same argument of the preceding proposition, and replace eigenvectors with modulus ``greater than 1" by ``not less than 1".  We  decompose the subspace $P_{2}$ into two subsubspaces $P_{2}''$ and $P_{2}'$, where $P_{2}''$ is generated by eigenvectors with the modulus of corresponding eigenvalues greater than 1, and $P_{2}'$ is generated by eigenvectors with the modulus of corresponding eigenvalues equal to 1. Then  $\R^{n}=P_{0}^{\pitchfork}\oplus P_{1}\oplus P_{2}'\oplus P_{2}''$, and let $\Pi_{2}':\R^{n}\rightarrow P_{2}'$ be the canonical projection induced by the decomposition above. So we also obtain the following inequality

\begin{center}
$
\begin{aligned}
& \parallel \Pi_{2}'(\tilde{f}^{n}(x)-\tilde{f}^{n}(y))\parallel \\
&\geq \parallel \Pi_{2}'(f_{\ast}(\tilde{f}^{n-1})(x)-f_{\ast}
(\tilde{f}^{n-1})(y))\parallel -\\
&\quad \parallel \Pi_{2}'(\tilde{f}^{n}(x)-\tilde{f}^{n}(y))-\Pi_{2}'
(f_{\ast}(\tilde{f}^{n-1})(x)-f_{\ast}(\tilde{f}^{n-1})(y))\parallel\\
&\geq \parallel \Pi_{2}'(\tilde{f}^{n-1}(x)-\tilde{f}^{n-1}(y))\parallel-c\\
&\ldots\\
&\geq  \parallel \Pi_{2}'(x-y)\parallel-nc.
\end{aligned}
$
\end{center}
Then we have

\begin{equation}
\parallel\tilde{f}^{n}(x)-\tilde{f}^{n}(y)\parallel \geq \parallel \Pi_{2}'(x-y)\parallel-nc. \label{(1)}
\end{equation}

Since $\tilde{\F}^{s}$ is quasi-isometric (see Propositon \ref{iso1}), for any $y\in \tilde{\F}^{s}(x)$ and each $n$ we can obtain that
\begin{equation}
 \parallel\tilde{f}^{n}(x)-\tilde{f}^{n}(y)\parallel \leq
k_{0}\lambda^{n}\cdot\parallel x-y\parallel, \label{(2)}
\end{equation}
where constant $k_{0}$ is independent of $x$ and $y$, and constant $\lambda<1$ is the contracting rate of $\tilde{f}$ on $\tilde{\F}^{s}$. Then, inequalities \eqref{(1)} and \eqref{(2)} imply that

\begin{equation}
\parallel \Pi_{2}'(x-y)\parallel-nc \leq k_{0}\lambda^{n}\cdot\parallel x-y\parallel,
 \label{(3)}
\end{equation}
for any $y\in \tilde{\F}^{s}(x)$ and each $n$.

Next, choose $y$ in $\tilde{\F}^{s}(x)$ such that $\parallel \Pi_{2}'(x-y)\parallel = tc$, where $t$ is a positive integer with $t>-2\log_{\lambda}2k_{0}k_{1}$, and  $k_{1}\triangleq\frac{\parallel x-y\parallel}{\parallel \Pi_{2}'(x-y)\parallel}$ (indeed, by Proposition \ref{bound}, there is a uniformly bounded distance between  $\tilde{\F}^{s}(x)$ and $P_{0}$, then we can always choose proper $y$ in $\tilde{\F}^{s}(x)$ such that the constant $k_{1}$ is bounded and $\parallel \Pi_{2}'(x-y)\parallel$ can be arbitrarily large). Then we take $n$ to be $\frac{t}{2}$ in inequality \eqref{(3)}, this inequality will be inverted, which is a contradiction. Hence, the modulus of eigenvalues corresponding to these eigenvectors generating $P_{0}$ are greater than one. Combining  Proposition \ref{esti}, it follows that $f_{\ast}$ is Anosov.

\end{proof}
\end{prop}

Since $f_{\ast}$ is Anosov, by a well-known result of J. Franks \cite{Fr}, there exists a semiconjugacy $H: \R^{n}\rightarrow \R^{n}$ which is $C^0$-close to $id$, and satisfies that for any $\alpha\in\Z^n$,  $H(x+\alpha)=H(x)+\alpha$ and
\begin{center}
$
H\circ \tilde{f}=f_{\ast}\circ H.
$
\end{center}
Since $H$ is $\Z^{n}$-periodic, $H$ can descend to $\T^{n}$, set it $h_{f}$.

Hence, $f$ is semiconjugate to $f_{\ast}$, the desired  codimension one Anosov diffeomorphism. The item $(iii)$ of Theorem \ref{final} holds.

Finally, we have completed the proof of Theorem \ref{final}.

\section{The dynamical structure}\label{sec5}

This section is devoted to the proof of Theorem \ref{dynachara}.  Let $f:\T^{n}\rightarrow \T^{n}$ be a codimension one partially hyperbolic diffeomorphism,  and let $H$ be the semiconjugacy between $\bar{f}$ and $f_{\ast}$ as mentioned at the end of the preceding section.

Since $f_{\ast}$ is Anosov, let $\kL^{s}$ and $\kL^{u}$
denote the stable and unstable foliations of $f_{\ast}$ in $\R^{n}$, respectively. Next, we will analyze the dynamical structure and behavior of this kind of system.

\begin{lem}\label{cusu}
The notations are as mentioned above, then

\begin{enumerate}[(1).]
   \item $H(\tilde{\F}^{u})\subseteq \kL^{u},\ i.e.,\ \forall x\in\R^n,\ H(\tilde{\F}^{u}(x))\subseteq \kL^u(H(x))$;
   \item $H(\tilde{\F}^{c})\subseteq \kL^{s},\ i.e.,\ \forall x\in\R^n,\ H(\tilde{\F}^{c}(x))\subseteq \kL^s(H(x))$.
\end{enumerate}

\begin{proof}
(1) Pick $y\in\tilde{\F}^{u}(x)$. Since $d_{\tilde{\F}^{u}}(\tilde{f}^{n}(y),\tilde{f}^{n}(x))\rightarrow \infty$ when $n\rightarrow\infty$, then so is $d(H(\tilde{f}^{n}(y)),H(\tilde{f}^{n}(x)))$ by $d(H, id)<const$,  which implies that $f_{\ast}^{n}(H(y))\in\kL^{u}(f_{\ast}^{n}(H(x)))=f_{\ast}^{n}(\kL^{u}(H(x)))$.
 Therefore, $H(\tilde{\F}^{u}(x))\subseteq \kL^{u}(H(x))$.

(2) Suppose the contrary: there is $x\in\R^{n}$ such that $H(\tilde{\F}^{c}(x))\nsubseteq\kL^{s}(H(x))$. Consider the iteration of $H(\tilde{\F}^{c}(x))$ over $f_{\ast}$, then $f_{\ast}^{n}(H(\tilde{\F}^{c}(x)))=H(\tilde{f}^{n}(\tilde{\F}^{c}(x)))=
H(\tilde{\F}^{c}(\tilde{f}^{n}(x)))$  will converge to $\kL^{u}(f_{\ast}^{n}(H(x)))$ in $C^{0}$-topology as $n$ increases, which is contradictory to $d(H,id)<const$ and $\tilde{\F}^{c}(\tilde{f}^{n}(x))$ transverse to $\tilde{\F}^{u}(\tilde{f}^{n}(x))$.

\end{proof}
\end{lem}

\begin{rmk}

The second relation implies that $\tilde{\F}^{c}$ can be uniformly bounded by a one dimensional linear foliation $\kL^{s}$ of $f_{\ast}$ in $\R^{n}$.

\end{rmk}

The above lemma allows us to show that $H$ can only collapse center arcs.

\begin{prop}\label{prei}
For any $z\in\R^{n}$, the preimage $H^{-1}(z)$ is a connected subset (i.e., an arc or a point) contained in a certain leaf of $\tilde{\F}^{c}$.
  \begin{proof}
Firstly, we prove that the preimage of each point must be contained in a leaf of $\tilde{\F}^{c}$. Take $z\in\R^{n}$ and pick $x,\ y\in H^{-1}(z)$. If $x\notin \tilde{\F}^{c}(y)$, then $H(\tilde{\F}^{c}(y))\subseteq \kL^{s}(H(y))=\kL^{s}(z)$ and $H(\tilde{\F}^{c}(x))\subseteq \kL^{s}(H(x))=\kL^{s}(z)$ by Lemma \ref{cusu}. It is easy to  check that $H|_{\tilde{\F}^{u}(x)}$ is injective. Indeed, for any $y,\ z\in \tilde{\F}^{u}(x)$, $H(y)=H(z)$ implies $y=z$. Otherwise, $f_{\ast}^{n}(H(y))=f_{\ast}^{n}(H(z))$, namely $H\circ \tilde{f}^{n}(y)=H\circ \tilde{f}^{n}(z)$ for each $n\geq 0$, which is impossible by $d(H,id)<const$.

\begin{figure}[htbp]
		\centering
		\includegraphics[width=11cm]{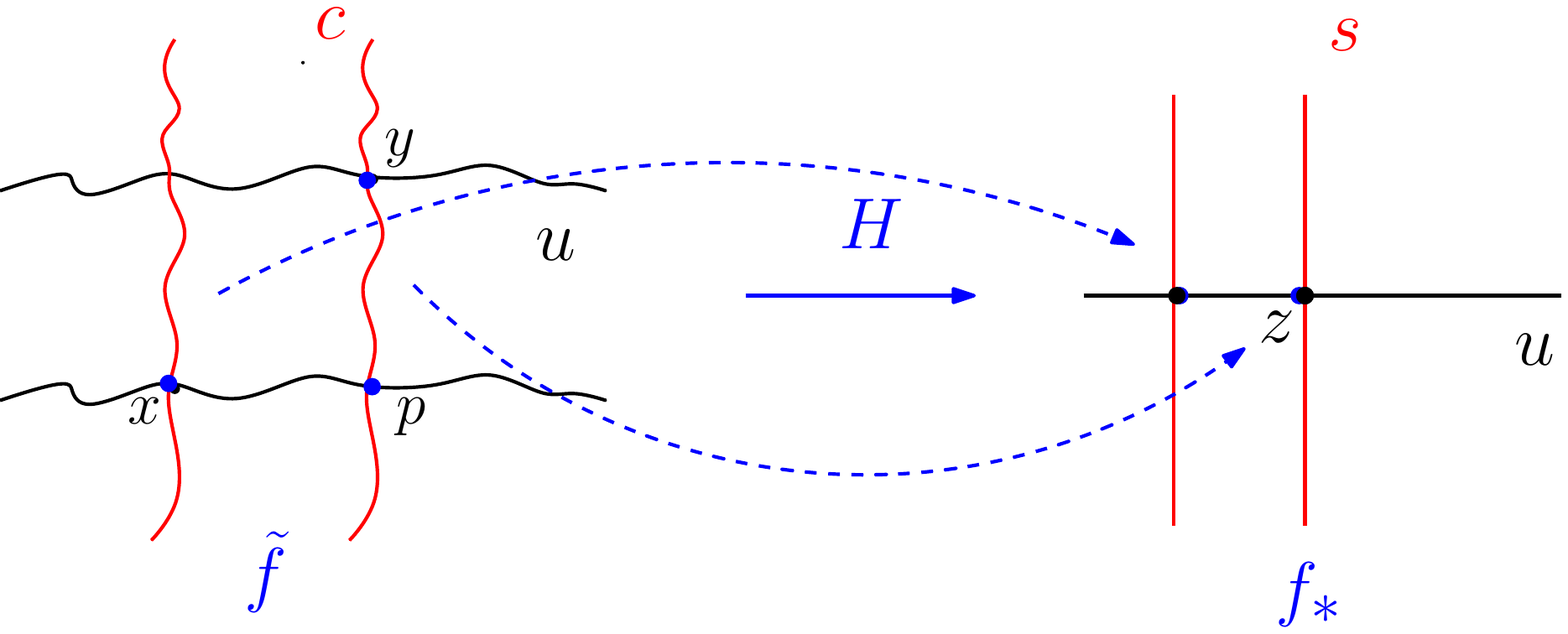}
		\caption{The diagram}
        \label{fig:graph}
\end{figure}

Set $p\in\tilde{\F}^{c}(y)\cap \tilde{\F}^{u}(x)$  as shown in  Figure \ref{fig:graph}, then $H(p)\in \kL^{s}(z)\setminus \{z\}$. When $n \rightarrow \infty$,  $d(\tilde{f}^{n}(x),\tilde{f}^{n}(p))\rightarrow \infty$, then $d(H\circ\tilde{f}^{n}(x),H\circ\tilde{f}^{n}(p))\rightarrow \infty$, but $d(f^{n}_{\ast}(z), f^{n}_{\ast}(H(p)))\rightarrow 0$, which is a contradiction.

For any $z\in \R^{n}$, $H^{-1}(z)$ is a compact set lying in a center leaf. Then we just need to prove the connectedness.  Observe that  for any $y,\ x\in \R^{n}$ and $H(y)=H(x)$ if and only if there exists a constant $C>0$ such that $d(\tilde{f}^{n}(y),\tilde{f}^{n}(x))<C$ for all $n\in\Z$, and $C$ can be taken to be independent of $y$ and $x$. Thus, when  $y$ and $x$ belong to $H^{-1}(z)$ we have $d(\tilde{f}^{n}(y),\tilde{f}^{n}(x))<C$. Then take $w$ in the center segment joining $y$ and $x$. Because $\tilde{\F}^{c}$ is quasi-isometric by Proposition \ref{ciso}, there are constants $a$ and $b$ such that for all $n\in \Z$, and we have

\begin{center}
$
\begin{aligned}
d(\tilde{f}^{n}(y),\tilde{f}^{n}(w)) & \leq d_{\tilde{\F}^{c}}(\tilde{f}^{n}(y),\tilde{f}^{n}(w)) \\
&\leq d_{\tilde{\F}^{c}}(\tilde{f}^{n}(y),\tilde{f}^{n}(x))\\
&\leq
ad(\tilde{f}^{n}(y),\tilde{f}^{n}(x))+b\\
&\leq aC+b.
\end{aligned}
$
\end{center}

Hence, we have $H(w)=H(y)=H(x)=z$, which follows that the whole center segment joining $y$ and $x$ is contained in $H^{-1}(z)$. So $H^{-1}(z)$ is connected.

  \end{proof}
\end{prop}

Next, we will apply the following lemma to prove Theorem \ref{dynachara}.

\begin{lem}[\cite{RP3}, Proposition 2.1]\label{wild}
If the semiconjugacy $H$ between $\tilde{f}$ and $f_{\ast}$ satisfies that:
\begin{enumerate}[(a).]
		\item $H$ is injective in unstable foliation $\tilde{\F}^{u}$;

        \item the fibers of $H$ are invariant under unstable holonomy maps;

        \item there exists a chain-recurrence $\mathcal{Q}$ in $\T^{n}$ such that the frontier of fibers of $H$ in center stable leaves are all contained in $\pi^{-1}(\mathcal{Q})$.

\end{enumerate}
Then every chain-recurrence class in $\T^{n}$ different from $\mathcal{Q}$ is contained in the preimage of a periodic orbit by $h_{f}$, the map of $H$ descending to $\T^{n}$.

\end{lem}

Notice that we have already proved that $H$ satisfies the above condition $(a)$, see the proof of Proposition \ref{prei}.  So we just need to demonstrate the following two lemmas, which show that $H$ also satisfies the above conditions $(b)$ and $(c)$.

\begin{lem}\label{invu}
 The fibers of $H$ are invariant under unstable holonomy maps.
\begin{proof}
Consider $x=H(y)$ for $y\in \R^{n}$. Since $d(H,id)<c$ and $\tilde{f}^{-n}(H^{-1}(x))=H^{-1}(f_{\ast}^{-n}(x))$, then for every $n>0$  we obtain that $diam(\tilde{f}^{-n}(H^{-1}(x)))<2c$.

It follows that there exists $n_{0}$ such that if $n>n_{0}$ then $\tilde{f}^{-n}(\Upsilon^{u}_{y,z}(H^{-1}(x)))$ is close enough to $\tilde{f}^{-n}(H^{-1}(x))$, where $z\in\tilde{\F}^{c}(y)$ and  $\Upsilon^{u}_{y,z}:\tilde{\F}^{c}(y)\rightarrow\tilde{\F}^{c}(z)$ denote the unstable holonomy map (because $\tilde{\F}^{c}$ and $\tilde{\F}^{u}$ have global product structure,  $\Upsilon^{u}_{y,z}$ can be  globally defined on the whole leaf $\tilde{\F}^{c}(y)$). In particular, we can choose a constant $\gamma>0$ such that if $z\in\tilde{\F}^{u}_{\gamma}(y)$ we have
\begin{center}
  $
diam(\tilde{f}^{-n}(\Upsilon^{u}_{y,z}(H^{-1}(x))))<4c.
  $
\end{center}

 To prove the lemma holds, it is sufficient to prove that $H(\Upsilon^{u}_{y,z}(H^{-1}(x)))$ contains exactly one point. Otherwise, the stable coordinates of $f_{\ast}$ of these points in $H(\Upsilon^{u}_{y,z}(H^{-1}(x)))$ must be different. Then, after backwards iterations, we will find that the distance between them can be much larger than $4c$. Since $H\circ \tilde{f}^{-n}=f_{\ast}^{-n}\circ H$ and $H$ is $c$-$C^{0}$-close to the identity $id$, there is a contradiction.

\end{proof}
\end{lem}

\begin{lem}\label{unia}
There is a unique quasi-attractor $\cQ$ for $f$. Moreover, every point $x$ which belongs to the frontier of a fiber of $H$ relative to its leaf $\tilde{\F}^{c}$ belongs to $\cQ$.
  \begin{proof}
Firstly, from the Introduction, we know that there exists a quasi-attractor for $f$.

Take $x\in\R^{n}$, and consider a point $y$ that lies on the boundary of $H^{-1}(\{x\})$ relative to $\tilde{\F}^{c}(y)$. Then because of lemma \ref{cusu} and the choice of $y$, we can know that the image of $y$ by $H$ cannot be contained in the unstable set of $x$ for $f_{\ast}$. Thus, iterating $\tilde{f}$ backwards, we can get a connected set of arbitrarily large diameter in the direction of the stable eigenline of $f_{*}$. Consider a quasi-attractor $\cQ'$ of $f$, for any $\varepsilon>0$ it follows that $\tilde{f}^{-m}\left(B_{\varepsilon}(y)\right)$ intersects $\pi^{-1}(\cQ')$ for sufficiently large $m$. This allows us to construct an $\varepsilon$-pseudo-orbit from $y$ to $\cQ'$. Then $y\in\cQ'$, because $\cQ'$ is a chain-recurrence class.

Since $\cQ'$ is chosen arbitrarily and all quasi-attractors are disjoint, we obtain that there is a unique quasi-attractor $\cQ$ for $f$.

  \end{proof}
\end{lem}

Let $\mathcal{Q}$ be the unique quasi-attractor obtained in lemma above. By the lemma \ref{wild}, we can know that every chain-recurrence class of $f$ different from $\mathcal{Q}$ is contained in a periodic interval. And the periodic interval  is the preimage of a periodic orbit by $h_{f}$, so it is an interval contained in the center foliation $\F^{c}$ of $f$.

Hence, this concludes the proof of Theorem \ref{dynachara}.

\section{Entropy}\label{sec6}

In this section, we show that Theorem \ref{cord} is true. Let $f:\T^{n}\rightarrow \T^{n}$ be a codimension one partially hyperbolic diffeomorphism, $H$ and $h_{f}$ be defined in the previous section.

Let $n \in \mathbb{N}$ and $\varepsilon>0$. A finite subset $E\subset \T^{n}$ is $(n, \delta)$-separated if for $\forall x, y \in E$ there exists $i\in\{0,1,\ldots, n-1\}$ such that $d(f^{i}(x),f^{i}(y))>\varepsilon$. Denote
\begin{center}
$\E_n(f, X, \varepsilon)\triangleq\text{max} \{\# E |\ E \subset X\subseteq \T^{n}$ is $(n, \varepsilon)$-separated $\}$
\end{center}
and
\begin{center}$
\E(f, X)\triangleq\lim \limits_{\varepsilon \rightarrow 0} \limsup\limits_{n \rightarrow \infty} \frac{1}{n} \ln \E_n(f, X, \varepsilon)$.
\end{center}
And when $X=\T^{n}$, $\E_{\text {top }}(f)\triangleq\E(f, \T^{n})$ is called the \emph{topological entropy} of $f$. We remark that $\E(f, X)=0$ if $X$ is a subset of a curve with all its iterates  having uniformly bounded length.

And in this section, we employ definition of entropy of an $f$-invariant measure $\rho$ (denote it as $\E_\rho(f)$) from ergodic theory consistent with \cite{Wal}. The well-known variational principle shows that $\sup \left\{\E_\rho(f) |\  \rho\right.$ is $f$-invariant$\}$ is exactly $\E_{\text {top }}(f)$. We say that an $f$-invariant measure $\rho$ is a \emph{ maximizing measure} if it satisfies $\E_\rho(f)=\E_{\text {top }}(f)$. For  general dynamical systems, the maximizing measure does not necessarily exist (see \cite{Mi}). Nevertheless, there always exists a uniquely maximizing measure on the basic sets of uniformly hyperbolic diffeomorphisms (see \cite{Bo}), and in the setting of this paper, a codimension one partially hyperbolic diffeomorphism $f$ always admits maximizing measures (see \cite{CY}\cite{DF}). Next, we will show that the maximizing measure of $f$ is unique  by applying some arguments in \cite{Ure}, which means that $f$ is intrinsically ergodic.

We need the  following lemma, which states that the injectivity of  $H$ holds almost everywhere on $\R^{n}$.

\begin{lem}\label{zero}
Let $\Sigma\triangleq\{x\in R^{n}|\ \# H^{-1}(x)>1\}$, then $m_{1}(\Sigma)=0$, 
where $m_{1}$ is the Lebesgue measure on $\R^{n}$. In particular, $m(\pi(\Sigma))=0$, $m$ is the Lebesgue measure on $\T^{n}$.
\begin{proof}
By Lemma \ref{cusu}, we know that $H^{-1}(\kL^{s}(z))$ is contained in $\tilde{\F}^{c}(x)$ for any $x\in H^{-1}(z)$, which implies that $H^{-1}(\kL^{s}(z))=\tilde{\F}^{c}(x)$. Set $\Sigma^{s}_{z}=\kL^{s}(z)\cap \Sigma$, then $H^{-1}(y)$ is a nontrivial interval of $\tilde{\F}^{c}(x)$ for all $y\in\Sigma^{s}_{z}$, which follows that $\Sigma^{s}_{z}$ is a countable set. Because of the arbitrariness of $z$ and the fact that $\kL^{s}$ is a one dimensional linear foliation, Fublini's Theorem implies that $m_{1}(\Sigma)=0$.

\end{proof}

\end{lem}

Now, we are able to prove  Theorem \ref{cord}.

\noindent\emph{The proof of  Theorem \ref{cord}} :
Firstly, we show that for any given $f_{\ast}$-invariant measure $\vartheta$ there is an $f$-invariant measure $\omega$ such that its $h_{f}$-image is $\vartheta$, $i.e.$, $\omega\circ h_{f}^{-1}=\vartheta$. Indeed, take a generic point $z$ for $\vartheta$ ($i.e.$, $\frac{1}{n} \sum_{i=0}^{n-1} \delta_{f_{\ast}^i(z)} \rightarrow \vartheta,\ \text {as} \ n \rightarrow \infty$), then pick $x\in h_{f}^{-1}(z)$, for $n\geq 1$ let
\begin{center}
$\displaystyle
\omega_{n}\triangleq\frac{1}{n} \sum_{i=0}^{n-1} \delta_{x}\circ f^{-i}.
$
\end{center}
Then, it is easy to know that

\begin{center}
$\displaystyle
\delta_{x}\circ h_{f}^{-1}=\delta_{z}\quad \text{and}\quad \omega_{n}\circ h_{f}^{-1}=\frac{1}{n} \sum_{i=0}^{n-1} \delta_{f_{\ast}^i(z)}.
$
\end{center}
Choose a convergent subsequence $\{\omega_{n_{j}}\}$ such that $\omega_{n_{j}}\rightarrow \omega$. Thus $\omega$ is an $f$-invariant measure with $\omega\circ h_{f}^{-1}=\vartheta$.

Notice that Lemma \ref{zero} implies the uniqueness of a measure $\mu$ whose $h_{f}$-image is the  Lebesgue measure $m$, the uniquely maximizing measure for $f_{\ast}$. So in order to prove $(i)$ we just need to demonstrate that the $h_{f}$-image of a maximizing measure of $f$ is $m$. Indeed, for any $f$-invariant measure $\rho$ and set $\hat{\rho}\triangleq \rho\circ h_{f}$, then $\hat{\rho}$ is an invariant measure for $f_{\ast}$. Since $f_{\ast}$ is a factor system of $f$, we have $\E_{\rho}(f)\geq\E_{\hat{\rho}}(f_{\ast})$. On the other hand, by the Ledrappier-Walters variational principle (see \cite{LW})
\begin{center}
$
\sup\limits_{\nu: \nu \circ h_{f}^{-1}=\hat{\rho}} \E_{\nu}(f)=\E_{\hat{\rho}}(f_{\ast})+\displaystyle\int_{\mathbb{T}^n} \E(f, h_{f}^{-1}(z)) d \hat{\rho}(z),
$
\end{center}
and $\E(f, h_{f}^{-1}(z))=0$ for each $z\in\T^{n}$ (because the $h_{f}$-preimages of points are arcs
of uniformly bounded length and the partition by $h_{f}$-preimages is $f$-invariant), then we obtain $\E_{\rho}(f)\leq\E_{\hat{\rho}}(f_{\ast})$. Thus,
\begin{center}
$
\E_{\rho}(f)=\E_{\hat{\rho}}(f_{\ast}).
$
\end{center}
It follows that $h_{f}$-image of the uniquely maximizing measure $\mu$ of $f$ is $m$. Obviously, $(f,\ \mu)$ and $(f_{\ast}, m)$ are isomorphic via $h_{f}$.

For $(ii)$, by \cite{KH}, one can know that the A. Katok's conjecture is true for Anosov diffeomorphism $f_{\ast}$. Since $h_{f}\circ f=f_{\ast}\circ h_{f}$ and the entropy of the fibers of $h_{f}$ is $0$ ($i.e.$, $\E(h_{f}^{-1}(z),\ f)=0$ for any $z\in\T^{n}$), as a consequence of the above Ledrappier-Walters variational principle, the A. Katok's conjecture  also holds for $f$ if for any given $f_{\ast}$-ergodic measure $\sigma$ there is an $f$-ergodic measure $\zeta$ such that $\zeta\circ h_{f}^{-1}=\sigma$. In fact, the previous statements imply that there is an $f$-invariant measure $\eta$ such that  its $h_{f}$-image is  $\sigma$. Next we just prove that the ergodic components of $\eta$ are desired $f$-ergodic measures. Observe that measure $\theta\circ h_{f}^{-1}$ is $f_{\ast}$-ergodic whenever measure $\theta$ is $f$-ergodic. Let $\eta=\int_{E(\mathbb{T}^{n},f)} \theta_{\tau} d \tau$ be the ergodic decomposition of $\mu$ (see Chapter 6 of  \cite{Wal}), where $E(\mathbb{T}^{n},f)$ denotes the set of all $f$-ergodic measures on $\mathbb{T}^{n}$. Acting the push-forward induced by $h_{f}$ to both sides of the above equation we can obtain that
\begin{equation}
\sigma=\eta\circ h_{f}^{-1}=\int_{E(\mathbb{T}^{n},f)} (\theta_{\tau}\circ h_{f}^{-1})d \tau. \label{(4)}
\end{equation}

According to the preceding observation, we get that for $\tau$-$a.e.$, $\theta_{\tau}\circ h_{f}^{-1}$ is $f_{\ast}$-ergodic measure. Hence, by the uniqueness of ergodic decomposition, \eqref{(4)}  implies that $\theta_{\tau}\circ h_{f}^{-1}=\sigma$, $\tau$-$a.e.$, where $\{\theta_{\tau}\}\subseteq E(\T^{n},f)$ are desired.

Henceforth, the proof of Theorem \ref{cord} is completed.

\hfill $\Box$\\

\acknowledgement  The author greatly appreciate the support of  Shaobo Gan and Yi Shi by numerous discussions, suggestions and valuable corrections. The author is very grateful to Jinhua Zhang for some advice, and Ruihao Gu, Yushan Jiang for helpful conversations with them. The author would also like to sincerely thank Rafael Potrie for kindly providing knowledge and references about codimension one foliations and the global product structure.

\bigskip


\end{document}